\documentclass[12pt]{article}
\usepackage{amsfonts}
\usepackage{epsfig}
\title{Square Turning Maps and their Compactifications}
\author{Richard Evan Schwartz \thanks{\hskip 5 pt Supported by 
N.S.F. Research Grant DMS-0072607}}

\newtheorem{theorem}{Theorem}[section]

\newtheorem{lemma}[theorem]{Lemma}

\newtheorem{corollary}[theorem]{Corollary}
\newtheorem{conjecture}[theorem]{Conjecture}

\def\startproof{{\bf {\medskip}{\noindent}Proof: }}

\def\endproof{$\spadesuit$  \newline}

\def\C{\mbox{\boldmath{$C$}}}% 
\def\Q{\mbox{\boldmath{$Q$}}}% 
\def\R{\mbox{\boldmath{$R$}}}% 
\def\T{\mbox{\boldmath{$T$}}}% 
\def\X{\mbox{\boldmath{$X$}}}% 
\def\Z{\mbox{\boldmath{$Z$}}}% 

\begin{document}
\maketitle
\begin{abstract}
In this paper we introduce some infinite
rectangle exchange transformations which are
based on the simultaneous turning of the squares
within a sequence of square grids.  We will show that such
noncompact systems have higher dimensional
dynamical compactifications. In good cases,
these compactifications 
are polytope exchange transformations
based on pairs of Euclidean lattices. In
each dimension $8m+4$ there is a
$4m+2$ dimensional family of them.  Here
$m=0,1,2,...$  The case $m=0$, which
we studied in depth in [{\bf S1\/}],
has close connections to the $E_4$ Weyl group and
the $(2,4,\infty)$ hyperbolic triangle group.
\end{abstract}

\section{Introduction}

\subsection{Background}

A {\it piecewise isometry\/} is a map defined on a
typically polyhedral subset of Euclidean space.  
The domain is partitioned into smaller
polyhedra in such a way that the map is defined,
and an isometry, when restricted to the interior
of each polyhedron in the partition.  
(One could make a similar definition
for piecewise affine maps.)  These maps often
have a special beauty and combinatorial feel to them,
but currently they are very far from being understood.
In particular, there is a scarcity of examples 
above dimension $1$, and especially above dimension $2$.

The simplest examples of piecewise isometries are $1$-dimensional
interval exchange transformations.
See [{\bf K\/}], [{\bf R\/}], [{\bf Y\/}], [{\bf Z\/}]
for an important but small sample of this large field.
The paper [{\bf H\/}] is an early paper on rectangle
exchange transformations.
The papers
[{\bf AE\/}],
[{\bf AG\/}],
[{\bf AKT\/}],
[{\bf Go\/}], 
[{\bf Low1\/}],
[{\bf Low2\/}], 
[{\bf LVK\/}],
[{\bf T\/}],
all treat closely
related sets of systems with $k$-fold
symmery for $k$ typically equal to
$5$, $7$, or $8$. 
Some definitive theoretical work concerning the
(zero) entropy of such maps is done in
[{\bf GH1\/}], [{\bf GH2\/}], and [{\bf B\/}].
My monographs [{\bf S1\/}] and
[{\bf S2\/}] have additional references.

The purpose of this paper is to study
the dynamics of maps defined in terms
of square grids. 
Given some $s>0$ and some $z \in \C$,
let $G_{s,z}$ be the (usual)
infinite grid of squares
in the plane such that $z$ is a vertex of
the grid and the sides of the squares are
parallel to $s$ and $is$. In particular, the
squares have side length $s$.
Let $R_{s,z}$ denote the map which rotates each
square of $G_{s,z}$ clockwise by $\pi/2$ radians.
The map $R_{s,z}$ is only defined on the interiors
of the squares of $G_{s,z}$. 

The map $R_{s,z}$ is not very interesting
from a dynamical perspective, because
$R_{s,z}^4$ is the identity map.  However,
the compositions of the form
\begin{equation}
\label{compose}
R_{S,Z}:=R_{s_1,z_1} \circ ... \circ R_{s_n,z_n}
\end{equation}
can have quite intricate behavior.
Here $S=\{s_1,...,s_n\}$ and
$Z=\{z_1,...,z_n\}$ are the data the the
composition.
When $n$ is divisible by $4$, the map
$R_{S,Z}$ is a piecewise translation.
More specifically, $R_{S,Z}$ is an infinite
rectangle exchange transformation in this case.

One nice feature of (planar) piecewise translations
is that they define a natural family of
convex polygons in the plane.
When $R=R_{S,Z}$ is a piecewise translation,
every periodic point $p$ of $R$
is contained in a maximal open convex
polygon $I_p$, called a {\it periodic island\/},
such that $R$ is defined entirely on $I_p$ and
periodic there.  Since $R$ is a rectangle
exchange, the polygon $I_p$ is in fact an
open rectangle.  

Probably the two most basic questions one can ask
about the maps in Equation \ref{compose} are
about their periodic points and about their
unbounded orbits.  In this paper we will talk
mainly about periodic points, though a few
of our results have to do with the existence
of unbounded orbits.  Our interest in understanding
these basic questions for a special case is
what led us to the idea of compactifying
the systems.

\subsection{A Motivating Example}

The original motivation behind this paper
was to understand the dynamics of the maps
\begin{equation}
A_s=(R_{1,0}R_{s,0})^2=R_{S,Z}, \hskip 30pt
S=\{1,s,1,s\}, \hskip 15 pt
Z=\{0,0,0,0\}.
\end{equation}
In [{\bf S1\/}], we called
$A_s$ the {\it alternating grid system\/}
because we generate the dynamics by alternately
turning the squares in one grid and in the other.

Let $Q_s$ denote the square of $G_{s,0}$ whose
bottom left vertex is the origin.
Figure 1.1 shows how $Q_s$
is tiled by the periodic islands for
$s=36/31$. The continued
fraction expansion of $s$ is $0:1:1:2:3:1:2$.

\begin{center}
\resizebox{!}{4.1in}{\includegraphics{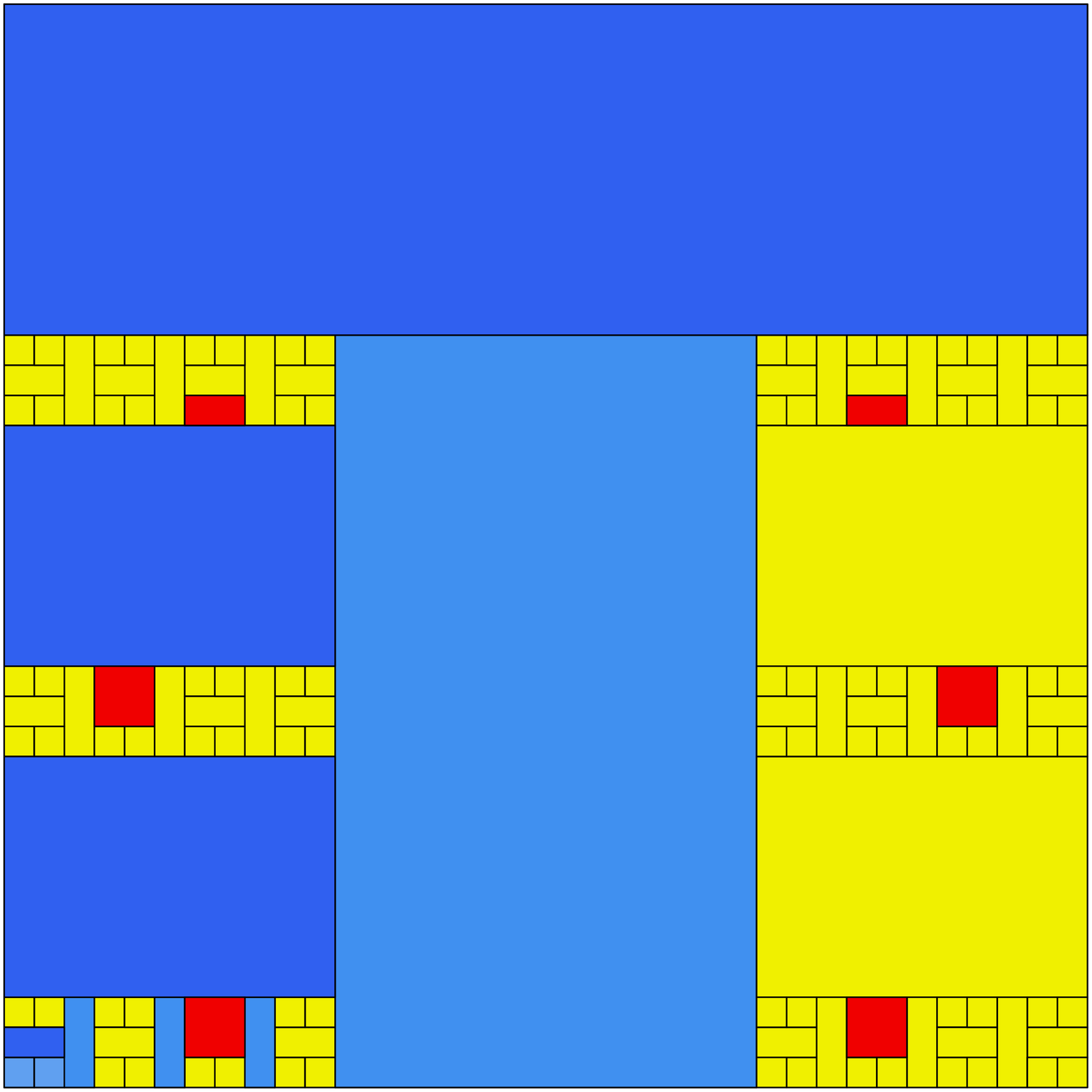}}
\newline
{\bf Figure 1.1:\/} Periodic island tiling of $Q_s$ for $s=36/61$.
\end{center}
 
If we study the blue islands in the picture,
we see that they encode the continued fraction
expansion of $s$.  The red islands suggest
``blemishes''.  It 
seems that these red islands should be subdivided
into smaller tiles in order to make the overall
pattern perfect.

Figure 1.2 shows the same kind of picture for the
parameter $s=55/89$.  Here, the C.F.E.
is $0:1:1:1:1:1:1:1:1:2$. 

\begin{center}
\resizebox{!}{4in}{\includegraphics{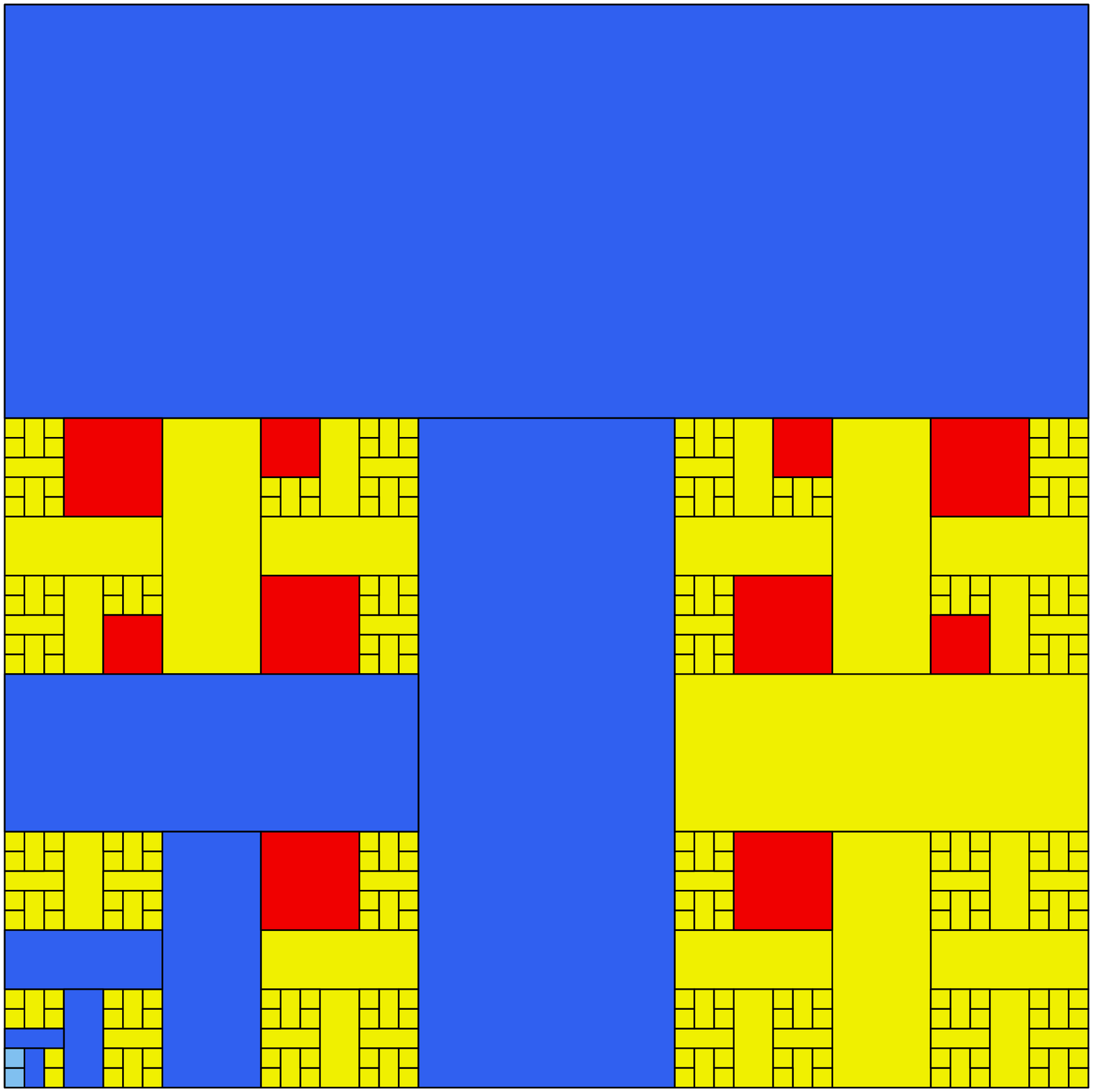}}
\newline
{\bf Figure 1.2:\/} Periodic island tiling of $Q_s$ for
$s=55/89$.
\end{center}

Pictures such as these suggest

\begin{conjecture}
\label{main}
Let $s \in R_+$.  Almost
every point of $A_s$ is periodic.
Inside $Q_s$ one can find a sequence of
periodic islands, alternately sharing
edges with the two coordinate axes,
which encodes the continued fraction
expansion of $s$.
\end{conjecture}

The difficulty in proving Conjecture \ref{main}
is that the periodic orbits corresponding to
small islands seem to be quite long and also to
have huge diameters. (The mysterious red
islands are also an obstacle to our understanding.)
This leads one to the idea
of trying to compactify the system, in order to
bring all the orbits into view, so to speak.
We took this approach (with partial success)
in [{\bf S1\/}] for $A_s$ and
here we will do it more generally.

\subsection{Dynamical Results}

The alternating grid system is an example of a
composition from Equation \ref{compose} which
just involves $2$ grids.  For the purpose
of stating some results, we formalize the idea
of such compositions.  Also, for convenience,
we always take the parameter $s$ to be irrational.
In the rational case, all the orbits are periodic.
There is a lot to say about the rational case,
but we will not say it here.

Let $G$ be the free product
$(\Z/4)*(\Z/4)$.  Let $A_1$ and $A_2$ be the
usual generators of $G$.  Let $G^0 \subset G$
denote the index $4$ subgroup consisting
of words whose total exponent is divisible
by $4$.  The word $(A_1A_2)^2$, which
corresponds to the alternating grid system,
is an example of a word in $G^0$.

Given $S=\{s_1,s_2\}$ and $Z=\{z_1,z_2\}$,
we have a representation of
$G$ into ${\rm map\/}(\C)$ which sends 
$A_j$ to $R_{s_j,z_j}$.  By scaling
and translating the plane, an operation
which does not effect the dynamics, it
suffices
 to consider the case when
$s_1=1$ and  \footnote{Somewhat later on, we will
find it more useful to set $z_1=(1+i)/2$, to that
the origin is the center of a square of the
first grid, rather than a vertex of the grid.}
$z_1=0$.  We set
$s=s_2$ and $z=z_2$.  We let 
$g_{s,z}$ denote the image of $g \in G$ under
our representation.  
When $g \in G^0$, the map $g_{s,z}$ is a piecewise
translation.

There is an auxilliary homomorphism
$H: G \to {\rm isom\/}(\C)$
such that $H(A_1)$ rotates by $\pi/2$ radians
counterclockwise around the origin and
$H(A_2)$ rotates by $\pi/2$ radians
counterclockwise about some other point.  We usually
take this other point to be $1$.
The homomorphism $H$ does not depend on the
parameters $s$ and $z$.   We say that the
word $g \in G^0$ is a {\it statonary word\/} if
$H(g)$ is the identity. We say that
$g \in G^0$ is a {\it drifter\/} if
$H(g)$ is not the identity map.
The set of stationary words is the kernel of $H$,
an infinite index subgroup of $G^0$.

There is a big difference between the stationary
words and the drifters.

\begin{theorem}
\label{fixed}
Suppose that $g \in G^0$ is a stationary word,
and $s$ is irrational, and $z$ is arbitrary. Then
$g_{s,z}$ has a positive density set of fixed
points.  In particular, $g_{s,z}$ has infinitely
many periodic islands of period $1$.
\end{theorem}

\begin{theorem}
\label{drift1}
Let $g \in G^0$ be a drifter and let $s$ be irrational.
For any $N>0$ and any $z$, the map
$g_{s,z}$ has only finitely many periodic islands
of period less than $n$.
Moreover, for all but countably many
choices of $z$, the map
$g_{s,z}$ has no periodic points at all.
\end{theorem}

The word $g=(A_1A_2)^2$ is a stationary word.
In this case, the results in [{\bf S1\/}] give
us an improved version of Theorem \ref{fixed}.

\begin{theorem}
\label{stationary1}
Let $g=(A_1A_2)^2$.
Let $s \in \R_+$ be irrational and
let $z \in \C$ be arbitrary.
There is an unbounded
sequence $\{p_k\}$ such that a positive
density set of points in $\C$
are periodic with period $p_k$ with
respect to $g_{s,z}$.
\end{theorem}
Theorem \ref{stationary1} takes a step in the
direction of Conjecture \ref{main}.
We will deduce Theorem \ref{stationary1}
from some results in [{\bf S1\/}], and we will give a
self-contained proof for the parameter $s=\sqrt 2/2$.

The same circle of techniques used to prove
Theorem \ref{stationary1} also proves

\begin{theorem}
\label{stationary2}
Let $g=(A_1A_2)^2$.
Let $s \in \R_+$ be irrational and
let $z \in \C$ be arbitrary.
Then for any $N$ there are orbits
having diameter greater than $N$.
Moreover, there exists uncountably
many choices of $z$ such that
$g_{s,z}$ has unbounded orbits.
\end{theorem}
Again, we will sketch the proof of Theorem
\ref{stationary2} in general, and give a self-contained
proof for the parameter $s=\sqrt 2/2$.

\subsection{The Compactifications}

Our proofs of the dynamical results mentioned
in the previous section rely on the
construction of compactifications for the
square turning systems.  Given $s \in \R_+$
and $z \in \C$, we get a group action of
$G=(\Z/4)*(\Z/4)$ on $\C$, as discussed in the
previous section. Our first result below gives a
compactification of this group action.

Let $\Z[i]$ be the Gaussian
integers.\footnote{We often
find it more convenient to work in
$\C^k$ rather than $\R^{2k}$, even
though sometimes we will make the
identification of $\C^k$ with
$\R^{2k}$.} 
Geometrically, $(\Z[i])^k$ is just the
square grid in $\C^k$. 
Let $\T^4=\C^2/(\Z[i]^2)$ be the usual 
$4$-dimensional square torus.
We define
$\Psi=\Psi_s: \C \to \T^4$ by the equation
\begin{equation}
\label{compact0}
\Psi(x+iy)=(z,z/s)\ {\rm mod\/}\ (\Z[i])^2.
\end{equation}
When $s$ is irrational $\Psi$ is injective on
$\C$ and $\Psi(\C)$ is dense in $\T^4$.
Note that $\Psi$ depends on $s$, but we
often suppress this dependence from our
notation.

\begin{theorem}[Compactification]
\label{compactify}
\label{one}
Suppose that $s$ is irrational.
There is a piecewise affine action
$\widehat G_s$ of $G$ on $\T^4$ having
the following properties.
\begin{enumerate}
\item For any $z \in \C$ there is some
translation $\widehat \tau$ of $\T^4$ so that
$\widehat \tau \circ \Psi$ conjugates the
action of $G_s$ on $\C$ to the
action of $\widehat G_s$ on
$\widehat \tau \circ \Psi(\C)$.
\item For each $g \in G$, the element
$\widehat g_s$ acts on
$\T^4$ in such a way that its
linear part $L_{g,s}$ 
independent of the point where
it is evaluated. 
\item When $g \in G^0$ is stationary,
$L_{g,s}$ is the identity, so that
$\widehat g_s$ is a polytope exchange
transformation.
\item 
When $g \in G^0$ is a drifter,
$L_{g,s}$ is a nontrivial complex linear
parabolic.  The fixed space of $L_{g,s}$ is
parallel to $\Psi(\C)$.
\end{enumerate}
\end{theorem}

Here is how we will deduce the dynamical results
from the Compactification Theorem
\begin{itemize}

\item When $g$ is stationary, the element
$\widehat g_s$ fixes a nontrivial polytope
in $\T^4$.  The set $\widehat \tau \circ \Psi(\C)$
intersects this polytope in a positive density set.
This is how we will prove
Theorem \ref{fixed}.  

\item When $g$ is a
drifter, the nontrivial parabolic nature of the action
of $\widehat g_s$ produces a local shear which
destroys the periodic points.
This is how we will prove
Theorem \ref{drift1}.

\item When $g=(A_1A_2)^2$ the map
$\widehat g_s$ has a compact invariant
$2$-dimensional slice, which we call an
{\it octagonal PET\/}.  In [{\bf S1\/}]
we showed that this invariant slice admits
a renormalization scheme.  This renormalization
scheme produces periodic points of arbitrarily
high order and also aperiodic points.  We
will deduce Theorems \ref{stationary1} and
\ref{stationary2} by studying how $\widehat \tau(\C)$
interacts with the orbits produced by the
renormalization scheme.

\end{itemize}

Independent of the dynamical consequences,
one might wonder about the geometry of the
compactifications produced by the
Compactification Theorem.  In some cases,
we can give a nice answer.
In \S 4 we will define what we mean by a
{\it double lattice PET\/}.  These are
special polytope exchange transformations
which are defined in terms of a pair of
Euclidean lattices and a pair of
fundamental domains for those lattices.
We first defined these maps in [{\bf S3\/}],
but [{\bf S1\/}] has a more thorough account.
The work in [{\bf S3\/}] is what led us to
the present paper.

\begin{theorem}
\label{double2}
Let $S=\{s_1,...,s_n\}$ and
$Z=\{z_1,...,z_n\}$, where
$n$ is congruent to $2$ mod $4$.
Let $R=R_{S,Z}^2$.
Then there is a double lattice
PET $(X,\widehat R)$, whose
domain $X \subset \R^{2n}$ is
a parallelotope, and an injective map
$\Psi: \C \to X$ which conjugates
$(R,\C)$ to the restriction of
$\widehat R$ to $\Psi(\C)$.
The closure of $\Psi(\C)$ is a finite
invariant 
union of $2k$-dimensional convex polytopes,
where $k$ is the dimension of the
$\Q$-vector space
$\Q[s_1,...,s_n]$.
\end{theorem}

The nicest case of our construction is when
there are no rational relations amongst $s_1,...,s_n$.
In this case, $\Psi(\C)$ is dense in $X$, and
$(X,\widehat R)$ is a compactification of
$(\C,R)$ in the most basic sense.  
Near the other
extreme, we can take periodic sequences
$S=\{1,s,1,s,...,1,s\}$ and
$Z=\{0,z,0,z,...,0,z\}$, where $s$ is irrational.
In this case, we get the following
corollary, which refers to the torus
action produced by the Compactification Theorem.

\begin{corollary}
\label{double}
Suppose that $g \in G_0$ has the form
$g=h^2$ where $h$ has word length 
$n=4m+2$.  Then the action of
$\widehat g_s$ on $\T^4$ is conjugate,
by a piecewise translation, to the
restriction of 
a $(2n)$-dimensional
double lattice PET to a finite
invariant union of $4$-dimensional
convex polytopes.
\end{corollary}

The double lattice PETs from Theorem \ref{double2}
occur in dimensions of the form $8m+4$ for $m=0,1,2,...$.
In dimension $8m+4$ there is a $4m+2$ dimensional family. 
\footnote{There is a certain redundancy, in the sense
that the data $\{\lambda s_1,...,\lambda s_n\}$ leads to
a system which is conjugate to the system defined
by $\{s_1,...,s_n\}$. So, perhaps it is more accurate
to say that there is a $4m+1$ dimensional family of
examples in dimension $8m+4$.}  We will describe our examples
explicitly for $m>0$ in \S 5, and for $m=0$ in \S 6.

For $m=0$, the examples have a special beauty; they
are related to the $E_4$-Weyl group and to the
$(2,4,\infty)$ hyperbolic triangle group.
 Our monograph
[{\bf S1\/}] is devoted to explaining the structure
associated to a $2$-dimensional invariant slice
of the compactification.  For
the purposes of explaining Theorems \ref{stationary1}
and \ref{stationary2}, we will re-derive some of the
structure here.
Here is a conjecture which encapsulates some of the
connection to the $E_4$ Weyl group.

\begin{conjecture}
\label{EE4}
When $g=(A_1A_2)^2$ and $s$ is
irrational, the associated double lattice PET
from Theorem \ref{double}
has an almost everywhere defined invariant
tiling $\widehat T_s$ by polytopes.  Each polytope 
in $\widehat T_s$ has order $192$ symmetry coming from
the action of the $E4$ Weyl group.
For any $z \in \C$ there is a leaf $L_z$ of the
invariant foliation such that 
$L_z \cap \widehat T_s$ is a refinement
of the tiling of $\C$ by periodic islands
of $g_s$.
\end{conjecture}
Referring to the discussion surrounding
Figures 1.1 and 1.2, the refinement 
we mention would be the result of suitably subdividing
all the ``red islands'' in the tiling, so as to improve
the picture.
Just so that this conjecture doesn't seem completely
off the wall, we will give a proof for the portions
of the tilings associated to fixed point sets.

\subsection{Organization}

This paper is organized as follows.
\begin{itemize}
\item In \S 2, we will prove 
the Compactification Theorem, Theorem \ref{one}.

\item  In \S 3 we prove Theorems \ref{fixed}
and \ref{drift1}.

\item In \S 4 we prove Theorem \ref{double}.

\item In \S 5 we will give an explicit description
of the PETs produced by Theorem \ref{double},
in dimensions $4,12,20,28,...$.

\item  In \S 6 we will consider in detail
the $4$-dimensional case from \S 5. 
We prove Theorems \ref{stationary1}
and \ref{stationary2} from
Theorem \ref{double} and 
a result
from [{\bf S1\/}],  
Theorem \ref{export}.
At the end of \S 6 we will prove the
special case of Conjecture \ref{EE4}
corresponding to the fixed point set.

\item In \S 7 we will roughly
sketch the proof of Theorem
\ref{export}, and we will give
a self-contained proof (modulo
a calculation) for $s=\sqrt 2/2$.

\item In \S 8 we discuss connections
to the $E_4$ Weyl group and prove the
special case of Conjecture \ref{EE4}
corresponding to the fixed point set.

\end{itemize}

\subsection{Acknowledgements}

I would like to thank Nicolas Bedaride,
Pat Hooper, Injee Jeong,
John Smillie, and
Sergei Tabachnikov for interesting
conversations about topics related to
this work.  I wrote this paper during
my sabbatical at Oxford in
2012-13.  I would especially like to
thank All Souls College, Oxford, for
providing a wonderful research environment.

My sabbatical was funded from many sources.
I would like to thank the National
Science Foundation, All Souls College,
the Oxford Maths Institute,
the Simons Foundation, the Leverhulme Trust, the
Chancellor's Professorship, and
Brown University for their support during this time
period.

\newpage

\section{The Compactification Theorem}

\subsection{The Main Construction}
\label{mainc}

Let $\T^{2n}$ denote the usual $2n$-dimensional torus
in $\C^n$.  We think of a fundamental domain for
$\T^{2n}$ as the unit cube centered at the origin
in $\C^n$.  We let $Q^{2n}$ denote the interior
of this cube.

Given $s$, define 
\begin{equation}
\label{goodchoice}
z=\frac{s+is}{2}.
\end{equation}
The grid with parameters $(s,z)$ is such that the origin
the center of one of the squares.
Let $R_s$ denote the corresponding grid map.

Given a finite sequence $s_1,...,s_n$, we set
$R_k=R_{s_k}$ and consider the composition
\begin{equation}
g=R_n \circ ... \circ R_1.
\end{equation}

We define $\Psi: \C \to \T^{2n}$ by the map
\begin{equation}
\Psi(z)=\bigg(\frac{z}{s_1},...,\frac{z}{s_n}\bigg)\ 
{\rm mod\/}\ (\Z[i])^n
\end{equation}
The image $\Psi(\C)$ is dense in a linear subspace
of $\T^{2n}$.  The dimension of this subspace, which we
will discuss in more detail below depends on the
number of rational relations between $s_1,...,s_n$.

We define
\begin{equation}
\label{R1}
\widehat R_1(z_1,...,z_n)=
(iz_1,z_2,...,z_n) + (-1+i) z_1
\bigg(0,\frac{s_1}{s_2},...,\frac{s_1}{s_n}\bigg)
\hskip 15 pt {\rm mod\/}\ (\Z[i])^n.
\end{equation}
This map is the identity on $\{0\} \times \C^{n-1}$
and locally affine on the set 
$Q^2 \times \T^{2n-2}$.
It is not possible to extend $\widehat R_1$ to
all of $\T^{2n}$, but this does not bother us.
The set $(\partial Q^2) \times \T^{2n-2}$, 
considered as a subset of $\T^{2n}$, is a
union of $2$ flat tori of dimension $2n-1$.
This is the singular set for $\widehat R_1$.

We define $\widehat R_j$ by permuting the coordinates and
interchanging the roles of $s_1$ with $s_j$.  For instance,
\begin{equation}
\label{R2}
\widehat R_2(z_1,...,z_n)=
(z_1,iz_2,...,z_n) + (-1+i) z_2
\bigg(\frac{s_2}{s_1},0,...,\frac{s_2}{s_n}\bigg)
\hskip 15 pt {\rm mod\/}\ (\Z[i])^n.
\end{equation}
The map $\widehat R_j$ is locally affine on the product
of an open square and a torus of dimension $2n-2$.

\begin{lemma}
$\widehat R_j \circ \Psi=\Psi \circ R_j$ on the
domain of $R_j$.
\end{lemma}

\startproof
By symmetry, it suffices to consider the case of $R_1$.
Notice that our lemma is true for the sequence
$\lambda s_1,...,\lambda s_n$ if and only if it
is true for the original sequence. For this reason,
it suffices to consider the case when $s_1=1$.
The domain of $R_1$ is the infinite grid of open unit
squares, one of which, namely $Q^2$, is centered at the origin.

We first check that the equation holds in a neighborhood
of $0$.  For $|z|$ sufficiently small, we compute
\begin{equation}
\Psi \circ R_1(z)=\Psi(iz)=
\bigg(iz,\frac{iz}{s_2},...,\frac{iz}{s_n}\bigg).
\end{equation}
On the other hand
$$
\widehat R_1(\Psi(z))=
\widehat R_1\bigg(z,\frac{z}{s_2},...,\frac{z}{s_n}\bigg)=$$
\begin{equation}
\bigg(iz,\frac{z}{s_2},...,\frac{z}{s_n}\bigg)+
(-1+i)z\bigg(0,\frac{1}{s_2},...,\frac{1}{s_n}\bigg)
\end{equation}
One can see that the two expressions are equal.

Next, we observe that $\Psi(Q^2) \subset Q^2 \times \T^{2n-2}$.
The restrictions of $R_1$ and $\Psi$ to $Q^2$ are affine
and the map $\widehat R_1$ is affine on
$Q^2 \times \T^{2n-2}$.  Therefore, the check we have
already made implies that our lemma holds true on $Q^2$.

Suppose we knew that $\Psi \circ R_1(z)=\widehat R_1 \circ \Psi(z)$
for some Gaussian integer $z$.  The differentials of 
$\Psi \circ R_1$ and $\widehat R_1 \circ \Psi$ at $z$ are
the same as the differentials of these maps at $0$.
Since they already agree at $0$, they agree at $z$
as well. This implies that our basic equation holds in
a neighborhood of $z$.  Note finally that
$\Psi$ maps the unit square $Q^2_z$ centered at $z$
into $Q^2 \times \T^{2n-1}$.  The same continuation
principle now implies that our lemma holds on
$Q^2_z$.
So, to finish the proof, we just have to check the
basic equation on Gaussian integers.

When $z$ is a Gaussian integer, we have $R_1(z)=z$.
So, we just have to prove that $\widehat R_1$ fixes $\Psi(z)$.
But $\Psi(z) \subset \{0\} \times \T^{2n-2}$, and such
points are fixed by $\widehat R^1$. 
\endproof

We define
\begin{equation}
\label{groupaction}
\widehat g=\widehat R_n \circ ... \circ \widehat R_1.
\end{equation}
An immediate consequence of our previous result is that
\begin{equation}
\label{conjX}
\Psi \circ g = \widehat g \circ \Psi
\end{equation}
wherever all maps are defined.

\subsection{Statement 1 of Theorem \ref{one}}

Theorem \ref{compactify} concerns the case $n=2$
in the construction above. In this case,
we normalize so that $s_1=1$ and $s_2=s$.
We will first prove Statement 1 for
the choice of $z$ in Equation \ref{goodchoice}.
In this case, Equation \ref{groupaction} gives
us the desired group action on $\T^4$.
When $z=0$, we take the translation $\widehat \tau$
in Theorem \ref{one} to be the identity.
Equation \ref{conjX} gives the desired conjugacy
between the two group actions.  This is Statement 1.

Now we explain how things work for other
choices of $z$.
Note that the grids corresponding to the
parameters $z$ and $z+sm+isn$ are the same,
for any $m,n \in \Z$.  Hence, it suffices
to prove Statement 1 for one representative
of each point in $\C$ in the torus
\begin{equation}
\T^2_s=\C/\Lambda, \hskip 30 pt \Lambda=\{sm+isn|\ m,n \in \Z\}.
\end{equation}
We first show that it suffices to consider a dense
set of points in the torus, and then we give an
argument which covers such a dense set of points.

\begin{lemma}
Suppose Statement 1 of the Compactification Theorem
holds for a dense set of points in $\T^2_s$.  Then
Statement 1 of the Compactification Theorem holds
for all points in $\C$.
\end{lemma}

\startproof
Let $z_{\infty} \in \T^2_s$ be some point, and let
$\{z_n\}$ be a sequence of points in our dense
set which converges to $z_{\infty}$.  Let $\widehat \tau_n$
be the corresponding set of translations of
$\T^4$.   Since the group of translations of
$\T^4$ is compact, we may pass to a subsequence
to that $\widehat \tau_n$ converges to a translation
$\widehat \tau$.   For any $n$, we have the equation
\begin{equation}
\label{equiv}
\Psi_n \circ G_{s,n}=\widehat G_s \circ \Psi_n,
\hskip 30 pt \Psi_n=\widehat \tau_n \circ \Psi.
\end{equation}
Here $G_{s,n}$ is the group action corresponding
to $z_n$.  Since $\Psi_n$ is an injection, 
Equation \ref{equiv} is equivalent to Statement 1 of
the Compactification Theorem for the parameter $z_n$.

Consider what happens as $n \to \infty$. For any
given element $g \in G$, the maps defined by
Equation \ref{equiv} converge uniformly on compact sets
to the corresponding maps defined in terms of the
limit parameter $z$.  Hence, Equation \ref{equiv}
holds as well for the $n=\infty$.  But then
Statement 1 of the Compactification Theorem
holds for $z_{\infty}$.
\endproof

Let $z_{00}$ be as in Equation \ref{goodchoice}.  Consider
the points
\begin{equation}
z_{mn}=z_{00}+m+in, \hskip 30 pt m,n \in \Z.
\end{equation}

What makes these points special is that the
group action defined in terms of $z_{mn}$ is
conjugate to the group action defined in
terms of $z_{00}$.  It is not generally true
that the group actions defined in terms of
different choices of $z$ are conjugate.

When $s$ is irrational, the set of representatives
of the points $\{z_{mn}\}$ in $\T^2_s$ is dense.  
So, it suffices
to prove Statement 1 of the Compactification Theorem
for the points $z_{mn}$.

\begin{lemma}
\label{explicit}
Statement 1 of the Compactification Theorem holds
for $z_{mn}$.
\end{lemma}

\startproof
Let $G_{s,m,n}$ be the group action corresponding
to the point $z_{mn}$.
Let $\tau_{mn}$ be the translation of $\C$ which
carries $z_{00}$ to $z_{mn}$.  By construction,
$\tau_{mn}$ conjugates $G_{s,0,0}$ to $G_{s,m,n}$.
More specifically,
\begin{equation}
\label{conj}
\tau \circ  G_{s,m,n} \circ \tau^{-1} = G_{s,0,0}.
\end{equation}
Since $\Psi$ is locally affine, there is some
translation $\widehat \tau_{mn}$ such that
\begin{equation}
\Psi \circ \tau_{mn}=\widehat \tau_{mn} \circ \Psi.
\end{equation}

For ease of notation we make
the following abbreviations.
\begin{equation}
\label{trans}
\tau=\tau_{mn}, \hskip 15 pt
\widehat \tau=\widehat \tau_{mn}, \hskip 15 pt
G=G_{s,0,0}, \hskip 15 pt
G'=G_{s,m,n}, \hskip 15 pt
\widehat G=\widehat G_s.
\end{equation}
Combining Equations \ref{conj} and \ref{trans} with
the special case of Statement 1 we have already
proved, we have
\begin{equation}
\widehat \tau \Psi G' \tau^{-1}=
\Psi \tau G'  \tau^{-1} =\Psi G=
\widehat G \Psi.
\end{equation}
Hence
\begin{equation}
\widehat \tau \Psi G'=
\widehat G \Psi \tau =
\widehat G \widehat \tau \Psi.
\end{equation}
This last equation says that
$\widehat \tau_{m,n} \circ \Psi$ is a
semi-conjugacy between $G_{s,m,n}$ and
$\widehat G_s$.  This equivalent to 
Statement 1 for $z_{mn}$.
\endproof

\subsection{The Rest of Theorem \ref{one}}

Statement 2 of Theorem \ref{one} is immediate.
The linear parts of the maps
$\widehat R_1$ and $\widehat R_2$ are independent
of the point where they are evaluated, and so
the same goes for any word in these generators.

Statements 3 and 4 of the Compactification Theorem
are purely statements of linear algebra.
 The linear parts of $\widehat R_1$ and
$\widehat R_2$ are given by

\begin{equation}
L_1=\left[\matrix{
i&0 \cr -1/s+i/s & 1}\right]
\hskip 30 pt
L_2=\left[\matrix{
1&-s+is\cr 0&i}\right]
\end{equation}
We introduce the matrix
\begin{equation}
E=\left[\matrix{s&0 \cr 1&1}\right].
\end{equation}
$E$ is a basis of eigenvectors of $L_1$.
Setting $M_k=E^{-1}L_kE$, we have
\begin{equation}
M_1=\left[\matrix{i&0\cr 0&1}\right]
\hskip 30pt
M_2=\left[\matrix{i&-1+i \cr 0 &1}\right].
\end{equation}
As long as $s \not = 0$, the matrix $E$ is
nonsingular.
Notice that $M_1$ and $M_2$ do not depend on $s$.
So, the conjugacy class of some word in $L_1$ and $L_2$
does not depend of $s$.

A calculation shows that $M_1$ and $M_2$ both preserve
$\Pi=\C \times \{-1\}$.
Moreover,
$M_1$ rotates $\Pi$ by $\pi/2$ counterclockwise
about the point $(0,-1)$ and $M_2$ does the same
thing about the point $(1,-1)$.
Therefore, the auxilliary 
homomorphism $H$ discussed
in connection with Theorem \ref{one} can be
interpreted as the map which
carries $A_k$ (a generator of
$G=\Z/4*\Z/4$) to $M_k|\Pi$ for $k=1,2$.
In short
\begin{equation}
H(A_k)=M_k|\Pi, \hskip 30 pt k=1,2.
\end{equation}

So, $g$ is a stationary word if and only if
the corresponding word in $M_1$ and $M_2$
acts trivially on $\Pi$.  More formally,
let $g \in G^0$.  Let
$\mu_g$ be the word in $M_1$ and $M_2$
corresponding to $g$.  By construction, 
\begin{equation}
E^{-1}L_g E = \mu_g.
\end{equation}
Here $L_g$ is the linear part of $\widehat g$.

Note that $g \in G_0$ forces $\mu_g$ to
be the identity on $\C \times \{0\}$.
Suppose $H(g)$ is trivial.  Then
$\mu_g$ is the identity on both
$\C \times \{0\}$ and $\C \times \{-1\}$.
Since $\mu_g$ is complex linear, this
forces $\mu_g$ to be the identity.
Hence $L_g$ is the identity as well.
This proves Statement 3.

We use the same notation for Statement 4.
If $H(g)$ is not the identity, then
the nontrivial action of $\mu_g$ on $\Pi$
forces $\mu_g$ to be a nontrivial parabolic.
Since $\mu_g$ is complex linear and acts
as the identity on $\C \times \{0\}$, it
must be the case that $\mu_g$ is a parabolic
of real rank $2$, when interpred as acting
on $\R^4$. The same goes for $L_g$,
which is conjugate to $\mu_g$ over $\C$ and
{\it a forteriori\/} over $\R$.  This proves
Statement 4.

\newpage

\section{Existence of Periodic Points}

\subsection{Good Partitions}

We fix some drifter $g \in G^0$ and consider
$\widehat g_s$ acting on $\T^4$.
Recall that the unit cube $Q^4$ is our
fundamental domain for $\T^4$. Recall also
that $\widehat g_s$ has an invariant
foliation $F=F_s$ that is parallel to $\Phi_s(\C)$.

We say that a
a {\it small polytope\/} is
a convex polytope $P \subset Q^4$
such that no face of $P$ contains an open
subset of a plane in $F$.
We also insist that
$g(P) \subset Q^4$.
In other words, neither $P$ nor $g(P)$ is allowed
to cross the boundary of $Q^4$.
Given a piecewise affine map $h: \T^4 \to \T^4$,
we say that a partition $\Pi$ of
$\T^4$ is {\it good\/} for $h$ if
$\Pi$ consists of small polytopes
and $h$ is defined and affine on
the interior of each one.

\begin{lemma}
\label{good}
There is a good partition for $\widehat g_s$.
\end{lemma}

\startproof
We suppress the dependence on $s$ from our argument.
From the construction in  \S \ref{mainc},
we see that
$\widehat R_1$ and $\widehat R_2$ are defined
on the complement of a finite union of flat
$3$-tori.  These $3$-tori are transverse
to our foliation.  By
removing additional $3$-dimensional flats, including
the boundary of the unit cube,
we can find a partition which is 
simultaneously good
for all the words $\widehat R_i^j$
where $i=1,2$ and $j=\pm 1$.

We produce a good partition for any word in
$\widehat R_1$ and $\widehat R_2$ by
induction on the length of the word.
Suppose that $g=\widehat R_1h$ and
$\Pi_h=\{H_1,...,H_n\}$ is a good partition
for $h$.  Let
$\{A_1,...,A_m\}$ be a good partition for
$\widehat R_1$.
We define
\begin{equation}
G_k=\bigcup_{i=1}^m h^{-1}(h(H_k) \cap A_i)
\hskip 30 pt
k=1,...,n.
\end{equation}
Then $G_k$ is a partition of $H_k$ into finitely
many small polytopes such that
$g$ is defined on each one.
If $B$ is some polytope in this partition,
then $g(B) \subset \widehat R_1(A_i)$ for some $i$.
Hence $g(B)$ is also small.
The union of all the polyhedra in $G_k$ as
$k$ ranges from $1$ to $n$, gives the
desired partition for $g$.
A similar construction works for 
$\widehat R_i^jh$ for the other
relevant choices of $i$ and $j$.
\endproof

For use in the next chapter, we record the following
corollary of our construction above.

\begin{corollary}
\label{refine}
For any $n$,
there is a good partition $\Pi_n$ for $\widehat g_s^n$.
Moreover we can choose these partitions so that
$\Pi_n$ is a refinement of $\Pi_m$ when $n>m$.
\end{corollary}

\subsection{Proof of Theorem \ref{drift1}}

Let $g \in G^0$ be a drifter.  We first show that
$\widehat g_{s,z}$ has only finitely many
periodic islands of period less than $N$.
Since there are only finitely many positive
integers less than $N$, it suffices to prove
instead that $\widehat g_{s,z}$ has only
finitely many fixed points of period $N$.
Replacing $\widehat g_{s,z}$ by 
$g_{s,z}^N$, it suffices to prove the
result when $N=1$. 

That is, we want to prove
that $g_{s,z}$ has only finitely
many fixed islands. By {\it fixed island\/},
we mean a periodic island corresponding to a
point of period $1$.  Any two fixed islands
are either identical or have disjoint interiors.
We will prove the result for 
$z=(s+is)/2$, as in Equation \ref{goodchoice}.
The general case has the same proof, except
that the map $\widehat \tau \circ \Phi$
is used in place of the map $\Phi$, for
some suitable translation $\widehat \tau$.
For ease of notation we set
$g=\widehat g_{s,z}$.

Let $\Pi$ be a good partition for $\widehat g$.
By the Compactification Theorem,
$\widehat g$ preserves $F$ and the restriction
of $\widehat g$ to a suitable leaf of 
the foliation $F$ is
conjugate to the action of $g$ on $\C$.

Say that a {\it good disk\/} is an open
polygon of the form
$D=X \cap P$, where $X$ is a $2$-plane parallel 
to the foliation $F$.  We say that $D$ is
{\it fixed\/} if $\widehat g$ fixes $D$ pointwise. 

\begin{lemma}
If $D$ is a good disk but not a fixed good disk,
then $\widehat g$ has no fixed points in $D$.
\end{lemma}

\startproof
Recall that $D$ is parallel to the
eigenspace of the linear part of $\widehat g$, and
the eigenvalues are all $1$.  Hence, the restriction
of $g$ to $D$ is a translation.  This means
that $D$ is either a fixed disk or $g$ fixes no points
of $D$ at all.
\endproof

Now we come to the key structural result.

\begin{lemma}
\label{finitefix}
Each small polytope of $\Pi$ contains at most one
fixed good disk.  Hence, the set of fixed points of
$\widehat g$ is contained in finitely many good disks.
\end{lemma}

\startproof
Suppose that $\widehat g$ fixes two good disks in a good
polytope $P$.  The restriction of $\widehat g$ to
the interior of $P$ is affine.  Let $A$ be the affine
map which extends $\widehat g|P$.  By construction,
$A$ is the identity on two parallel $2$-planes.
But then $A$ is the identity on a certain $3$-dimensional
subspace.  But this is impossible, because the linear
part of $A$ is a parabolic of real rank $2$.
\endproof

Let $X$ denote the set of fixed points of $g$.
Let $X' \subset X$ denote these points $x$ such
that $\Phi(x)$ lies in the interior of some
small polytope of the good partition.  Since
the faces of these small polytopes are transverse
to the invariant foliation, we see that $X'$ is
dense in $X$.  In particular, every fixed
island contains points of $X'$.

Let $Y$ denote the set of fixed good disks.  Since
$\Psi_s: \C \to \T^4$ is injective and since all the
fixed points of $\widehat g$ lie in fixed good disks,
$\Psi_s$ induces a map from $X'$ into $Y$.

\begin{lemma}
Suppose that $p,q \in X'$.
If $\Psi$ maps $p$ and $q$ into the same fixed good
disk, then $p$ and $q$ lie in the same fixed island.
\end{lemma}

\startproof
Let $I_p$ and $I_q$ be the fixed islands
containing $p$ and $q$ respectively.
Either $I_p=I_q$ or these sets have disjoint interiors.
If $\Psi(p)$ and $\Psi(q)$ lie in the same good fixed disk $D$,
then $\Psi$ maps the line segment $\overline{pq}$ into $D$.
But then $\widehat g$ fixes every point of
$\Psi(\overline{pq})$.  Hence $\widehat g$ is entirely defined,
and the identity, on $\overline{pq}$.  
But then $\overline{pq} \subset I_p \cap I_q$.
This forces $I_p=I_q$.
\endproof

If $g$ had infinitely many fixed islands, then
the previous result would give us infinitely many fixed
good disks.   This is a contradiction.  Hence
$g$ only has finitely many fixed islands.
As we mentioned above, this completes the proof
that, for any $N$, the map $g=g_{s,z}$ has
only finitely many periodic islands of period less than $N$.
This is the first statement of Theorem \ref{drift1}.

Now we prove Statement 2 of Theorem \ref{drift1},
which says that $g_{s,z}$ has periodic
points only for countably many choices of $z$.
Using the same trick as for Statement 1, it
suffices to prove this result for fixed points.

Let $\widehat \tau_z$ be the translation of
$\T^4$ such that $\widehat \tau \circ \Psi_s$
conjugates the action of $g_{s,z}$ on $\C$
to the action of $\widehat g_s$ on 
\begin{equation}
\C_z=\widehat \tau_z \circ \Psi_s(\C).
\end{equation}
Here $\C_z \subset \T^4$ is one of the
leaves of the invariant foliation.
Note that the map $\C \to \C_z$ is a
bijection.  Here is the key lemma for our
result.

\begin{lemma}
If $g_{s,z}$ has a fixed point, then
$\C_z$ is one of finitely many leaves of the
invariant foliation $F$.
\end{lemma}

\startproof
There are only finitely many fixed good disks.
Hence there are only finitely many leaves of $F$
which contain good fixed disks.  Call these leaves
{\it special\/}.
If $g_{s,z}$ has a fixed point then $\C_z$ contains
an fixed good disk and hence is special.
\endproof

Again, we call a leaf in the invariant foliation
special if it contains a fixed good disk.
The next result shows
that there are only countable many choices of
$w$ such that $\C_w$
is a special leaf.  This finishes the proof
that there are only countably many choices
of $z$ for which $g_{s,z}$ has a fixed point.

\begin{lemma}
\label{count}
For each $z \in \C$, there are only countably
many choices of points $w \in \C$ such that
$\C_z=\C_w$.
\end{lemma}

\startproof
Suppose that $\C_z=\C_w$.
We can interpret the translation by
$\widehat \tau_z -\widehat \tau_w$ as
a translation of $\C$ which conjugates the action of
of $G_{s,z}$ to the action of $G_{s,w}$.
Call this translation $\tau$.

Since $\tau$ conjugates the rotation $R_1$ to the rotation
$R_1$, we see that $\tau$ must be translation by some
integer lattice vector.  In particular, $\tau$ preserves
the grid $G_1$.  At the same time, 
$\tau(G_{s,z})=G_{s,w}$.  Hence, there are only countably
many choices for $G_{s,w}$.  But only countably many
choices of $w$ lead to the same choice of $G_{s,w}$.
Hence, there are only countable many choices of $w$ which
lead to one of the countably many possible choices of grid.
\endproof

\subsection{Proof of Theorem \ref{fixed}}

Now suppose that $g \in G^0$ is a
stationary word.
Inspecting the maps $\widehat R_1$ and
$\widehat R_2$ which define the group
action $\widehat G_s$, we see that
both elements fix the origin $O_4 \in \T^4$
and both elements are defined in a neighborhood
of $O_4$. Hence, for any $g \in G$,
the word $\widehat g_s$
fixes $O_4$ and is defined in a
neighborhood of $O_4$.
When $g \in G^0$ is a stationary word, $\widehat g_s$ is a
piecewise translation.  In this case, there is
some nontrivial and maximal open polytope $P$
which contains $O_4$ in its interior such that
$\widehat g_s$ fixes every point of $P$.

We choose $z \in \C$ and let
$\Psi=\tau \circ \Psi_s$ be the map from the
Compactification Theorem.
Let $P'=\Psi^{-1}(P)$.  Since $\Psi$ is
locally affine and $\Psi(\C)$ is dense
in $\T^4$, we see that
$P'$ has positive density.  In fact, the
density of $P'$ is just the volume of $P$.
But $g_{s,z}$ fixes every point of $P'$.
Hence $g_{s,z}$ has a positive density
set of fixed points. Since $P'$ has
positive density, $P'$ is unbounded. Since
every fixed island is compact, there must
be infinitely many fixed islands.
This completes the proof of Theorem \ref{fixed}.

\subsection{A Generalization}

It seems worth mentioning a generalization of 
some of the results in this chapter.
Let $\widehat g$ be a piecewise affine map of
the flat torus 
$\T^n$.  We say that an {\it invariant foliation\/}
is a foliation $F$ such that $\widehat g$ preserves
every leaf of $F$.  We also reqire that 
$\widehat g$ is defined almost everywhere on
every leaf of $F$.  

\begin{theorem}
\label{shear}
Suppose $\widehat g$ is a piecewise affine map of $\T^n$
having a flat, dense, invariant $k$-dimensional
foliation $F$. Suppose that the linear part $L$ of $\widehat g$
is independent of the point of evaluation and is a parabolic
of real rank $k$ whose real eigenspace is tangent to $F$.
Then there are only finitely many leaves of $F$ containing
fixed points of $\widehat g$ and only countably many
leaves containing periodic points.
\end{theorem}

The proof of Theorem \ref{shear} is almost
the same as what we have done for
Theorem \ref{drift1}.
Mainly, we will point out the differences.
First of all, the first statement in
Theorem \ref{shear} can be applied to 
powers of $\widehat g$. Hence, the
first statement of Theorem \ref{shear}
implies the second statement.  Thus,
it suffices to prove that only finitely
many leaves of $F$ contain fixed points of
$\widehat g$.

The same argument as
in Lemma \ref{good} shows that $\T^n$ has a
good partition.  Here we use the fact that
$\widehat g$ is almost everywhere defined on
every leaf of $F$.  This guarantees that the
faces of the polyhedra in the initial
partition for $\widehat g$ do not contain
open subsets of the leaves of $F$.

There are two kinds of fixed points of
$\widehat g$, those contained in the interiors
of the small polytopes of the good partition,
and those contained in the faces of these
polytopes.  Call these fixed points of the
first kind and second kind respectively.

\begin{lemma}
Any fixed point of the second kind is contained
in the same leaf as some fixed point of the first kind.
\end{lemma}

\startproof
Let $p \in \T^n$ be a fixed point of the second kind.
If $p$ lies in a face of the good partition, it means
that $\widehat g$ is actually defined on $p$.  Since
$\widehat p$ is a parabolic of real rank $k$ whose
real eigenspace is parallel to $F$, we see that
$\widehat g$ fixes an entire $k$-disk $\Delta$
containing $p$.  By assumption, $\Delta$ is not
contained in the union of the boundaries of the small
polytopes.  Hence $\Delta$ contains fixed points of
the first kind.  But $\Delta$ lies in a single leaf
of $F$.
\endproof

In light of the previous result, we just have to
prove that the fixed points of the first kind are
contained in finitely many leaves of $F$.
Lemma \ref{finitefix} now goes through, almost word for
word, to show that $\widehat g$ has only finitely
many fixed good disks, and that these fixed good
disks contain all the fixed points of the first kind.
These finitely many fixed good disks lie inside
finitely many leaves of $F$.

\newpage
\section{Nature of the Compactifications}

\subsection{Double Lattice PETs}
\label{dlp}

In this chapter we prove Theorem
\ref{double2}.  For starters, we
define what we mean by a double lattice PET.

The data for a double lattice PET is a
quadruple $(X_1,X_2,\Lambda_1,\Lambda_2)$, where
\begin{itemize}
\item $X_1$ and $X_2$ are polytopes in $\R^n$.
\item $\Lambda_1$ and $\Lambda_2$ are lattices in $\R^n$.
\item $X_i$ is a fundamental domain for $\Lambda_j$
for all possible $i,j \in \{1,2\}$.
\end{itemize}
To say that $\Lambda$ is a lattice is to say that
there is some affine isomorphism $T$ of $\R^n$ such that
$\Lambda=T(\Z^n)$.  To say that $X$ is a fundamental
domain for $\Lambda$ is to say that the orbit
$\Lambda(X)$ tiles $\R^n$:  The translates have pairwise
disjoint interiors and the union of the translates
is a covering.

In all our examples, the polytopes $X_1$ and $X_2$
are parallelotopes centered at the origin.

There is a natural map $f_{ij}: X_i \to X_j$ defined
as follows.  Given $p \in X_i$, there is generically
a unique vector $\lambda_p \in \Lambda_j$ such that
$p+\lambda_p \in \Lambda_j$.  We define
\begin{equation}
f_{ij}(p)=p+\lambda_p.
\end{equation}
The maps
$f_{11}$ and $f_{22}$ are the identity; we do not 
care about these maps.  The maps $f_{12}$ and $f_{21}$
are the ones of interest to us.  These maps are piecewise
translations.  The composition
\begin{equation}
f=f_{21} \circ f_{12}: X_1 \to X_1
\end{equation}
is a polytope exchange transformation (PET) having
$X_1$ as a domain.  We call the system
$(X_1,f)$ a {\it double lattice PET\/}.  

We needed to break symmetry in order to define the maps
$f_{ij}$.  We could equally well define the maps $g_{ij}: X_i \to X_j$
as follows.  For $p \in X_i$ there is generically a unique
vector $\mu_p \in \Lambda_i$ such that $p+\mu_i \in X_j$.
We can then define $g_{ij}(p)=p+\mu_p$.
It is easy to see that
$g_{ij}=f_{ji}^{-1}$.  Hence
$f^{-1}=g_{12}g_{21}$.

At first it might seem difficult to produce quadruples
satiasfying the necessary conditions.  However, in
[{\bf S3\/}] we showed that essentially all polygonal
outer billiards systems have compactifications which
are double lattice PETs.  In this chapter we will
see that the construction in the previous
chapter leads naturally to double lattice PETs as well.

\subsection{The Linear Part}

We continue with the notation from \S \ref{mainc}
As a start on the proof of Theorem \ref{double2},
we analyze the linear part $L_g$ of $\widehat g$ in
case $g$ is a word of length $n=4m+2$.  Here
is our main result.

\begin{lemma}
\label{involution}
$L_g$ is an involution 
whose $(-1)$-eigenspace is $2$-dimensional
and whose $(+1)$-eigenspace is $(2n-2)$-dimensional.
\end{lemma}

We prove Lemma \ref{involution} through
a series of smaller results.  We also note
that the calculations done in \S \ref{dlp2} give
an alternate proof of Lemma \ref{involution}.

\begin{lemma}
The $(-1)$ eigenspace of $L_g$ has
real dimension at least $2$.
\end{lemma}

\startproof
The action of $\widehat g$ on $\Psi(\C)$ is conjugate
to the action of $g$ on $\C$.  Since $g$ has
length $4n_2$, the linear part of $g$ is
rotation by $\pi$. Hence, the same goes for $\widehat g$.
This implies that $L_g$ preserves the complex line
through the origin and parallel to $\Psi(\C)$.
In other words, considered as a real matrix
$L_g$ has a $(-1)$-eigenspace which is at least
$2$ dimensional.
\endproof

To finish the proof of
Lemma \ref{involution}, 
we will produce a $(2n-2)$-dimensional
subspace on which $L_g$ is the identity.  This will
finish the proof.
We will first consider the case when the numbers
$1,s_2,...,s_n$ have no rational relations
amongst them.  Once we take care of this case,
we will deduce the general case by a limiting argument.
When there are no rational
relations amongst the numbers $1,s_2,...,s_n$,
the image $\Psi(\C)$ is dense in $\T^{2n}$.

We say that a subset $\Sigma \subset \C$ is a
{\it net\/} if there is some $K$ such that every
point of $\C$ is within $K$ of some point of
$\Sigma$.  
If $U$ is any open subset of $\T^{2n}$, the
inverse image $\Psi^{-1}(U)$ is a net in $\C$.

\begin{lemma}
Let $\epsilon>0$ be given.  If $U$ is sufficiently small
then every point of $\Psi^{-1}(U)$ is within 
$\epsilon$ of a fixed point of $g$.
\end{lemma}

\startproof
We will crucually use the fact that the linear
part of $g$ is rotation by $\pi$.
Consider a point $z \in \C$ such that
$\Psi(z)$ lies within $\delta$ of the origin.
Then $z$ is very nearly the center of all the
grids used to define the maps
$R_1,...,R_n$.  Were $z$ the center of all
these grids, $g$ would locally be a rotation
by $\pi$ about $z$.  As it is, $g$ is a small
perturbation of a rotation by $\pi$ about $z$,
and the differential $dg$ is still rotation
by $\pi$. But, in this situation, $g$ has a
fixed point $z'$ very close to $z$.  
The distance $|z-z'|$ only depends on
$\delta$ and $n$.
\endproof

\begin{corollary}
Let $U$ be any open set containing the origin
in $\C^n$.  The set of fixed points of $g$
mapping into $U$ is a net in $\C$.
\end{corollary}

We choose some small open set $U$ as in the
corollary, and let $H$ denote the smallest
linear subspace containing all the points of
the set $\Psi(\Sigma)$.  By construction,
$\widehat g$ fixes every point of
$\Psi(\Sigma)$.  Since these points span $H$,
we see that $\widehat g$ is the identity on
$H$.  We just have to prove that $H$ has
dimension $2n-2$.

\begin{lemma}
$H$ has dimension at least $2n-2$.
\end{lemma}

\startproof
Suppose that $H$ has dimension $d<2n-2$.
Given any point $p \in \Sigma$, there is
some small disk $\Delta$ so that
the restriction of $g$ to $\Delta$ is
rotation by $\pi$.  The radius of $\Delta$
can be chosen to be uniformly large.
Call this radius $\rho$.
The image $\Psi(\Delta)$ is an isometric
disk centered at a point of $H$ and 
parallel to $\Psi(\C)$. 

Let $\Sigma'$ denote the $\rho$-tubular
neighborhood of $\Sigma$.  By construction,
$\Psi(\Sigma')$ is a union of disks, all
parallel to $\Psi(\C)$, and all centered
at points of $H$.  But this implies that
the closure of $K$ of
$\Psi(\Sigma')$ has dimension at most $d+2<2n$.

On the other hand, since $\Sigma$ is a net,
$\Sigma'$ is a set of positive density in
$\C$.  Given the affine nature of $\Psi$
and the fact that $\Psi(\C)$ is dense
in $\T^n$, we see that $K$ must have
positive volume -- i.e. $(2n)$-dimensional
Lebesgue measure.  This is impossible
if $\dim(K)<2n$.   This contradiction
shows that $\dim(H) \geq 2n-2$.
\endproof

We have produced a subspace of dimension at
least $2n-2$ which is fixed by 
$\widehat g$.  This means that the linear
part $L_g$ has a $(+1)$-eigenspace of
dimension at least $2n-2$.  This is all
we needed for Lemma \ref{involution} in
the arational case.

In case there are some rational relations
between the numbers $1,s_2,..,s_n$, we can
perturb the numbers slightly to get a new
sequence with no rational relations.
Hence, in the general case, the linear
part of $\widehat g$ has a sequence of
approximations by linear maps which all
have $(+1)$-eigenspaces of dimension
$2n-2$.  But then the dimension of the
$(+1)$-eigenspace of the limiting map
must be at least $2n-2$.  This completes
the proof of Lemma \ref{involution} in
the general case.

We mention the obvious corollary of
Lemma \ref{involution}

\begin{corollary}
Suppose that $g=h^2$ where $h$ is a word
of length $4m+2$.  Then the linear part
of $\widehat g$ is the identity. Hence
$(\T^{2n},\widehat g)$ is a PET.
\end{corollary}

\subsection{A Picture}

Before we continue our analysis in the
next section,
we show a picture that the reader should
keep in mind throughout the discussion.

The square shown in Figure 4.1 is meant to
be the torus $\T^2$.  The opposite sides are
meant to be identified.  The thick line
$\Omega_1$ is really a circle and the
thick line $\Omega_2$ is really a line
segment.  The set $\Omega_1 \cup \Omega_2$
is meant to be a toy version of the
singular sets for the map $\widehat h$ 
we consider in the next section.
The complement of the singular set is a
parallelogram -- i.e. a set affinely
eqiovalent to $Q^2$ -- embedded in $\T^2$
in a funny way.

\begin{center}
\resizebox{!}{4in}{\includegraphics{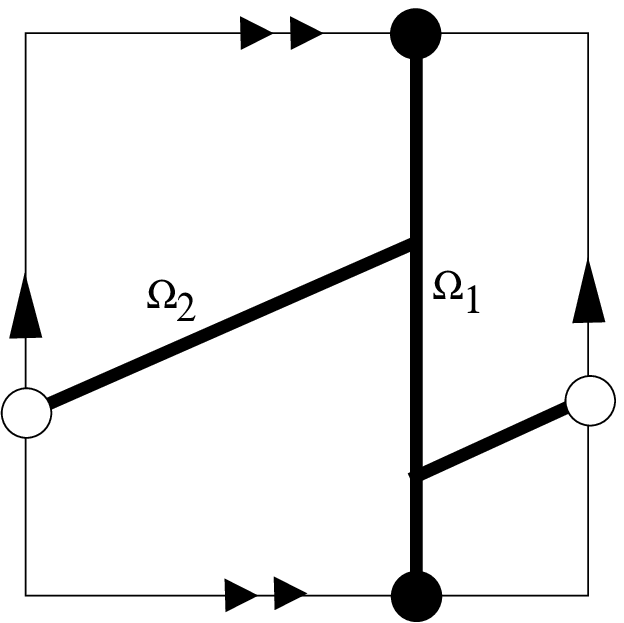}}
\newline
{\bf Figure 4.1:\/} A parallelogram sitting
inside a torus.
\end{center}

The set $\Omega_1$ is a $1$-torus sitting inside
$\T^2$.  The set $\Omega_2$ is obtained from
a horizontal line by applying a locally affine
map defined in the annulus $D_1=\T^2-\Omega_1$.
This locally affine map does not extend to
$\Omega_1$; it is a kind of partial Dehn twist.
The universal
cover $\Delta_1$ of $D_1$ is an infinite strip.
This strip is a low dimensional version of the
slabs $\Delta_k$ we consider in our proof below.

\subsection{The Singular Set}
\label{singularproof}

We continue with the case where $g=h^2$
and $h$ has length $n=4m+2$. 
Let $\Sigma_j \subset \T^{2n}$ denote the singular set
of $\widehat R_j$.  Recall that
$\Sigma_j$ is a union of $2$ 
codimension $1$ tori. 
Let $\Omega_1=\Sigma_1$ and then
\begin{equation}
\label{recursive}
\Omega_k=\widehat R_1^{-1} \circ ...
\circ \widehat R_{k-1}^{-1}(\Sigma_k), \hskip 30 pt
k=2,...,n.
\end{equation}
It follows from induction that
$\widehat h$ is entirely defined on
\begin{equation}
\T^{2n}-\Omega, \hskip 30 pt
\Omega=\bigcup_{j=1}^n \Omega_j.
\end{equation}

Recall that $Q^k$ is the open unit $k$-dimensional
cube centered at the origin.
The goal of this section is to prove the
following result.  

\begin{lemma}
\label{singular}
$\T^{2n}-\Omega$ is affinely equivalent to 
$Q^{2n}$.
\end{lemma}
We will prove Lemma \ref{singular} through
a series of smaller results.

Say that a {\it singular hyperplane\/} is a
hyperplane parallel to some $\Omega_j$.
Each $\Omega_j$ supplies $2$ singular hyperplanes,
so there are $2n$ singular hyperplanes in total.

\begin{lemma}
\label{transversality}
The $2n$ singular hyperplanes are linearly
independent in the sense that their normal
vectors form a real basis for $\R^{2n}$.
\end{lemma}

\startproof
We use complex notation.
Let $e_1,...,e_6$ denote the standard basis
vectors for $\C^6$. 
Let $L_j$ denote the linear part of
$\widehat R_j$. The hyperplanes associated
to $\Omega_1$ are
\begin{equation}
(e_1)^{\perp}, \hskip 5 pt
(ie_1)^{\perp}, \hskip 5 pt
L_1^{-1}((e_2)^{\perp}), \hskip 5 pt
L_1^{-1}((ie_2)^{\perp}), \hskip 5 pt
L_1^{-1}L_2^{-1}((e_3)^{\perp}), \hskip 5 pt
L_1^{-1}L_2^{-1}((ie_3)^{\perp}),...
\end{equation}
The corresponding normals are given by
\begin{equation}
\label{normal1}
e_1, \hskip 8 pt
ie_1, \hskip 8 pt
L_1^{t}(e_2), \hskip 8 pt
L_1^{t}(ie_2), \hskip 8 pt
L_1^{t}L_2^{t}(e_3), \hskip 8 pt
L_1^{t}L_2^{t}(ie_3),...
\end{equation}
An easy calculation shows that
the $k$th column of
$L_1^t...L_{k-1}^t$ has a $1$ in the
$(kk)$th position and $(0)$s below.
This implies that the normals in
Equation \ref{normal1} are linearly independent.
\endproof

We will prove Lemma \ref{singular} by
induction. Define
\begin{equation}
D_k=\T^{2n}-(\Omega_1 \cup ... \cup \Omega_k).
\end{equation}
Our final goal is to show that
$D_n$ is affinely equivalent to
$Q^{2n}$.  Our induction step will be that $D_k$
is affinely equivalent to
$\Q^{2k} \times \C^{n-k}$.
We have already seen that this statement is
true for $k=1$.
We will suppose that the statement is true for
some choice of $k$ and then prove it for
$k+1$.  

Define
\begin{equation}
\widehat h_k=\widehat R_k \circ ... \circ \widehat R_1.
\end{equation}
The map $\widehat h_k$ is defined on $D_k$.
We have projection
\begin{equation}
\pi: \C^n \to \T^{2n}.
\end{equation}
Let $L_k$ be the linear part of
$\widehat h_k$.
Let $\Delta_k \subset \C^n$ denote the universal cover of
$D_k$.   We think of $\Delta_k$ as being one connected
component of the preimage $\pi^{-1}(D_k)$.
The set $\Delta_k$ is affinely equivalent to the
``slab'' $Q^{2k} \times \C^{n-k}$.   

We have a commuting square
\begin{equation}
\matrix{\Delta_k& \stackrel{L_k}{\to} &\C^n 
\cr \downarrow \pi &&\downarrow \pi \cr 
D_k&\stackrel{\widehat h_k}{\to}&\T^{2n}}
\end{equation}

The set $\pi^{-1}(\Sigma_{k+1})$ consists of
two infinite families of parallel hyperplanes.
We call these infinite families
$A_{k+1}$ and $B_{k+1}$.  The two families are
transverse to each other.   By Lemma \ref{transversality},
the hyperplanes in
$L_k^{-1}(A_{k+1})$ and $L_k^{-1}(B_{k+1})$ are transverse to
$\partial \Delta_k$.

The intersection $\Delta_k \cap L_k^{-1}(A_{k+1})$ is an infinite
parallel family of sets, each of which is
affinely equivalent 
to $Q^{2k} \times \C^{n-k-1}$.  The same goes for
$D_k \cap L_k^{-1}(B_{k+1})$.  The projection $\pi$
carries these sets to $\Omega_{k+1}$.

Let $\Delta_{k+1}$ be some connected component of
$$\Delta_k-L_k^{-1}(A \cup B).$$
The set $\Delta_{k+1}$ is affinely
equivalent to 
$\Q^{2k+2} \times \C^{2n-2k-2}$.
Informally, the sets $L_k^{-1}(A_{k+1})$ and
$L_k^{-1}(B_{k+1})$ chop up two of the noncompact directions
into compact pieces.
The map $\pi: \Delta_{k+1} \to D_{k+1}$
is the universal covering map.  From this picture
we see that
$D_{k+1}$ is affinely equivalent to 
$\Q^{2k+2} \times \T^{2n-2k-2}$.  This
completes the induction step.

\subsection{The Double Lattice PET}

Now we start enhancing the notation from the previous
section.  Define
\begin{equation}
\Lambda_1=(\C[i])^n.
\end{equation}
Let $\pi: \C^n \to \C^n/\Lambda_1$ be projection, as above.

Since $\T^{2n}-\Omega$ is affinely equivalent
to the contractible set $Q^{2n}$, we have a
parallelotope $X_1$ and a continuous local inverse
\begin{equation}
\label{defX1}
\pi^{-1}: \T^{2n}-\Omega \to X_1 \subset \C^n.
\end{equation}
The parallelotope $X_1$ is just the image
of $\T^{2n}-\Omega$ under the lift
$\pi^{-1}$.
Since $\pi$ is injective on $X_1$ and
$\pi(X_1)$ has full measure in
$\T^{2n}$, we see that
$X_1$ is a fundamental domain for
$\Lambda_1$.

\begin{lemma}
\label{KEY}
$\widehat h=\pi \circ L_h \circ \pi^{-1}$
on $\T^{2n}-\Omega$.
\end{lemma}

\startproof
Both maps agree in a neighborhood of the origin and
are locally affine and entirely defined on the
contractible domain in question. Hence, these two
maps agree everywhere.
\endproof

Now we introduce another lattice and another
parallelotope.  Recall that the linear
part $L_h$ of $\widehat h$ is an affine
involution. We define

\begin{equation}
X_2=L_h(X_1), \hskip 30 pt
\Lambda_2=L_h(\Lambda_1).
\end{equation}
Now we have specified all the data for a double
lattice PET.  We just have to check the conditions.

\begin{lemma}
$X_2$ is a fundamental domain for $\Lambda_1$.
\end{lemma}

\startproof
We already know that
$\widehat h=\pi \circ L_h \circ \pi^{-1}$ on
$T^{2n}-\Omega$.  But this means that the
map
$\pi \circ L_h: X_1 \to \T^{2n}$
is injective and has dense image.
But this means that
$\pi: L_h(X_1) \to \T^{2n}$
is injective and has dense image.
Hence $L_h(X_1)$ is a fundamental
domain for $\Lambda_1$.
But $L_h(X_1)=X_2$.
\endproof

Since $L_h$ is an involution, we see that
$X_1$ is also a fundamental domain for
$\Lambda_2$. In short, the data
$(X_1,X_1,\Lambda_1,\Lambda_2)$ define a
double lattice PET.  We denote this
PET by $(X,f)$.  Here $X=X_1$ is the domain.

\begin{lemma}
The locally affine isomorphism
$\pi: X_1 \to \T^{2n}-\Omega$ conjugates the
system $(X,f)$ to the system $(\T^{2n},\widehat g)$.
\end{lemma}

\startproof
We first want to understand the map
$$
F=\pi^{-1} \circ \pi: X_2 \to X_1.
$$
For any $p \in X_2$, the points
$p$ and $F(p)$ differ by an
element of $\Lambda_1$.
Thus, the map $F_1: X_2 \to X_1$ is just the
map $f_{21}$ discussed in \S \ref{dlp}.
In short,
\begin{equation}
f_{21}=\pi^{-1} \circ \pi.
\end{equation}
Since $L_h$ is an involution which
interchanges the roles of
$X_1$ and $X_2$, and also the roles of
$\Lambda_1$ and $\Lambda_2$, we have
\begin{equation}
f_{12}=L_h \circ \pi^{-1} \circ \pi \circ L_h.
\end{equation}

Now we put these equations together.
\begin{eqnarray}
f= \cr \cr
f_{21} \circ f_{12}= \cr \cr
\pi^{-1} \circ \pi \circ L_h \circ \pi^{-1} \circ \pi \circ L_h=_1 \cr \cr
\pi^{-1} \circ \pi \circ L_h \circ \pi^{-1} \circ \pi \circ L_h
\circ (\pi^{-1} \circ \pi) = \cr \cr
\pi^{-1} \circ (\pi \circ L_h \circ \pi^{-1}) \circ (\pi \circ L_h
\circ (\pi^{-1}) \circ \pi =_2 \cr \cr
\pi^{-1} \circ \widehat h \circ \widehat h \circ \pi= \cr \cr
\pi^{-1} \circ \widehat g \circ \pi.
\end{eqnarray}
Equality 1 comes from the fact that
$\pi^{-1} \circ \pi$ is the identity on $X=X_1$.
Equality 2 is two applications of Lemma \ref{KEY}.
In short,
$f=\pi^{-1} \circ \widehat g \circ \pi$.
\endproof

It is worth emphasizing that the domain for $\pi^{-1}$ is
$\T^{2n}-\Omega$, though we think of
$\pi$ as giving a piecewise isometric conjugacy
between a system in $\T^{2n}$ and a system in
$X_1$.  Since our maps are not everywhere defined,
the difference in topology has no meaning here.

\subsection{The Invariant Slice}

Consider the composition
\begin{equation}
\Xi=\pi^{-1} \circ \Psi: \C \to X_1.
\end{equation}
$\Xi$ is not defined on the set
$\Psi^{-1}(\Omega)$.  This set is a
countable collection of line segments.
The fact that $\Xi$ is not defined at
these points does not bother us.  These
are the points where the original
planar system $g_{s,z}$ is not defined.
Putting together the results above,
we see that
$\Xi$ conjugates the action of $g_{s,z}$
on $\C$ to the action of the
double lattice PET $(X,f)$ on $\Xi(\C)$.

Now let $k$ denote the dimension of the
$\Q$-vector space $\Q[s_1,...,s_n]$.
Looking at the map $\Psi$, we see that
$\Psi(\C)$ is a $(2k)$-dimensional linear
subspace of $\T^{2m}$.  The intersection
\begin{equation}
\widehat C={\rm closure\/}(\Psi(\C)) \cap (\T^{2n}-\Omega)
\end{equation}
is a finite invariant union of $(2k)$-dimensional
open polytopes.
The set $\pi^{-1}(\widehat C)$ is likewise
an invariant finite union of convex
$(2k)$-dimensional polytopes in $X_1$.
At the same time, this set is the closure
of $\Xi(\C)$ in $X_1$.  This proves
Theorem \ref{double2}.

\subsection{Proof of Corollary \ref{double}}

It only remains to reconcile our compactification
here with the one from Theorem \ref{compactify}.
recall that $n=4m+2$.
In case we have the sequence
$1,s,1,s,...,1,s$, the construction here is 
identical to the construction given in
\S 2, except that we are repeating the coordinates
$2m+1$ times. In other words, the diagonal embedding
$\C^2 \to \C^m$ carries the compactification
produced by Theorem \ref{compactify} to the
system $(\widehat T^{2n},\widehat g)$.  The
composition of $\pi^{-1}$ with the diagonal
embedding gives the desired conjugacy.

\newpage

\section{A Concrete Family of PETs}

\subsection{Generalities}

We will work in $\C^n$.
We can specify \footnote{Not all double lattice
PETs can be specified this way, but the ones
here can be.} a double lattice PET
$(X_1,X_2,\Lambda_1,\Lambda_2)$ by a
quadruple of $n \times n$ matrices
$(\chi_1,\chi_2,L_1,L_2)$, where
\begin{itemize}
\item $X_j=\chi_j(Q^{2n})$.
\item $\Lambda_j$ is the $\Z[i]$-span
of the columns of $L_j$.
\end{itemize}
Here $Q^{2n}=(Q^2)^n$, where $Q^2$ is the unit
square centered at the origin in $\C$.

There is some ambiguity in our choice of 
matrices.
Let us call two vectors $V$ and $V'$ equivalent
if either $V'=\omega V$ or
$V'=\omega \overline V$.  Here
$\omega$ is some $4$th root of unity and
$\overline V$ is the complex conjugate of $V$.
Typically there are $8$ vectors in each
equivalence class.  We say that two matrices
$M$ and $M'$ are equivalent if
each column of $M$ is equivalent to the corresponding
column of $M'$.  If we replace the matrix
$\chi_j$ by an equivalent matrix $\chi_j'$
then we still recover $X_j$.  Likewise,
if we replace the matrix
$L_j$ by $L_j'$ we still recover $\Lambda_j$.

Here is a criterion which will help us verify that
the matrices we list give rise to double lattice PETs.
\begin{lemma}
Suppose that $\chi$ and $L$ are $n \times n$ matrices.
Let $X=\chi(Q^{2n})$ and let $\Lambda$ be the $\Z[i]$-span
of the columns of $L$.  Then $X$ is a fundamental domain
for $\Lambda$ provided that $\chi^{-1}L$ is a
triangular and has $4$-th roots of unity
along the diagonal.
\end{lemma}

\startproof
$X$ is a fundamental domain for $\Lambda$ if and only if
$T(X)$ is a fundamental domain for $T(\Lambda)$.  Here
$T$ is any invertible linear transformation of 
$\R^{2n} \approx \C^n$.  In particular, this is true
for $T=\chi^{-1}$.  In other words, it suffices to consider
the case when $\chi$ is the identity matrix and
$X=Q^{2n}$.  Replacing $L$ by an equivalent
matrix, we can assume that $L$ has $1$s along the
diagonal. But $Q^{2n}$ is indeed a fundamental domain
for a lattice whose defining matrix is
triangular and has $1$s along the diagonal.
\endproof

We call $(\chi_1,\chi_2,L_1,L_2)$ a {\it special system\/}
if $\chi_i^{-1}L_j$  is a triangular matrix with
$4$th roots along the diagonal for each pair $(i,j) \in \{1,2\}$.
The following corollary produces an $n$-parameter family of
double lattice PETs from a special system.

\begin{corollary}
Let $D$ be any nonsingular diagonal matrix. If
$(\chi_1,\chi_2,L_1,L_2)$ is a special system, then
$(\chi_1D,\chi_2D,L_1D,L_2D)$ is the data for a double lattice PET.
\end{corollary}

\startproof
We compute
\begin{equation}
(\chi_i D)^{-1}(L_jD)=D^{-1}(\chi_i^{-1}L_j)D=D^{-1}TD=T'.
\end{equation}
Here $T$ is the triangular matrix guaranteed by the
hypotheses, and $T'$ is the conjugate
matrix.  $T'$ is also triangular, and has $4$th roots
of unity along the diagonals.
So, for all indices $i,j$, the
parallelotope $\chi_i D(Q^{2n})$ is a fundamental domain for
the lattice defined by $L_jD$.
\endproof

\noindent
{\bf Remark:\/}
We could produce an example of a special
system by taking $4$ upper triangular matrices
(or lower triangular matrices).  However, this
would lead to a fairly trivial double lattice PET.
The map would essentially be a parabolic linear
transformation.

\subsection{An Explicit Example}
\label{exp000}

We discovered the construction in this section by
working out the details of the compactifications
described in the previous section.  We will present
the construction first and then identify it with what
we did in the previous chapter.  Our construction
works for $n=2,6,10,14,...$.
We introduce matrices
\begin{equation}
\chi_{\pm}=\left[\matrix{
\pm i&\pm i&\pm i&\pm i&\dots \cr
+i&-1& 0&0&\dots\cr
0&+i&-1& 0&\dots\cr
0&0&+i&-1&\dots\cr
\dots&\dots&\dots&\dots&\dots
}\right]
\end{equation}

\begin{equation}
L_{\pm}=\left[\matrix{
\mp 1&\pm i&\pm 1&\mp i&\dots \cr
+1&-1& 0&0&\dots\cr
0&+1&-1& 0&\dots\cr
0&0&+1&-1&\dots\cr
\dots&\dots&\dots&\dots&\dots
}\right]
\end{equation}

The first row for $L_{\pm}$ repeats every $4$ entries.
Each entry in this first row is $(-i)$ times the
preceding element.  Just to be clear, the top left
entry of $L_+$ is $(-1)$.
We set
$\chi_1=\chi_+$ and
$\chi_2=\chi_-$ and
$L_1=L_+$ and $L_2=L_-$.

Now we will verify that $(\chi_1,\chi_2,L_1,L_2)$ is
a special system.  Letting $J$ be the diagonal
matrix with entries $(-1,+1,+1,+1,...)$ we see
that $\chi_2=I\chi_1$ and $L_2=I L_1$. We compute
\begin{equation}
\chi_2^{-1}L_1=(I\chi_1)^{-1}L_2=\chi_1^{-1}I^{-1}L_1=
\chi_1^{-1}I L_1=\chi_1^{-1}L_2.
\end{equation}
Similarly,
\begin{equation}
\chi_2^{-1}L_2=\chi_1^{-1}L_1.
\end{equation}
For this reason, we just have to check the conditions
for $\chi_1^{-1}L_1$ and $\chi_1^{-1}L_2$.
A routine calculation verifies that
\begin{equation}
L_1=\chi_1\left[\matrix{
1&0&0&0&0&0&\dots\cr
-1+i&1&0&0&0&0&\dots\cr
-1-i&-1+i&1&0&0&0&\dots\cr
+1-i&-1-i&-1+i&1&0&0&\dots\cr
+1+i&+1-i&-1-i&-1+i&1&0&\dots\cr
-1+i&+1+i&+1-i&-1-i&-1+i&1&\dots\cr
\dots&\dots&\dots&\dots&\dots&\dots&\dots}
\right]
\end{equation}
The matrix listed is such that entry in
the lower triangle is $(-i)$ times the entry
directly below it.  Going down a
column, the pattern has period $4$.
A similar calculation verifies that
\begin{equation}
L_2=\chi_1\left[\matrix{
-i&-1+i&+1+i&+1-i&-1-i&-1+i&\dots\cr
0&-i&-1+i&+1+i&+1-i&-1-i&\dots\cr
0&0&-i&-1+i&+1+i&+1-i&\dots\cr
0&0&0&-i&-1+i&+1+i&\dots\cr
0&0&0&0&-i&-1+i&\dots\cr
0&0&0&0&0&-i&\dots\cr
\dots&\dots&\dots&\dots&\dots&\dots&\dots}
\right]
\end{equation}
The matrix listed is such that each entry in the
upper triangle is $(i)$ times the entry immediately
to its right.  Going across a row, the pattern 
has period $4$.

Thus, $(\chi_1,\chi_2,L_1,L_2)$ is a special system.
if $S=\{s_1,...,s_n\}$, we let
$D_S$ be the diagonal matrix whose diagonal
entries are $s_1,...,s_n$.
We will see that the double lattice PET
specified by
$(\chi_1D_S,\chi_2D_S,L_1D_S,L_2D_S)$ is 
conjugate to the one produced by Theorem
\ref{double2}.

\subsection{The Invariant Foliation}
\label{fol000}

Choose some $D_S$.  Let $X_j=\chi_jD_S(Q^{2n}))$ and
let $\Lambda_j$ be the lattice which is the
$\Z[i]$ span of the columns of $L_jD_S$.  
Here we discuss some symmetries of the
double lattice PET $(X_1,X_2,\Lambda_1,\Lambda_2)$.
Let $J$ be the involution mentioned in the
previous section.  By construction, the
action of $J$ swaps $X_1$ and $X_2$, and
simultaneously swaps $\Lambda_1$ and $\Lambda_2$.

Considering $J$ as a matrix acting on $\R^{2n}$,
this map has a $2$-dimensional $(-1)$-eigenspace
and an $(2n-2)$-dimensional $(+1)$-eigenspace.
This is just like the map $L_g$ considered in
connection with Theorem \ref{double2}.
The $(-1)$-eigenspace of $J$ defines a
complex line foliation $\cal F$ of $X_1$ (and $X_2$).
The leaves of this foliation are parallel to
$\C \times 0^{n-1}$.

\begin{lemma}
The foliation $\cal F$ is invariant under the
action of the double lattice PET.
\end{lemma}

\startproof
Let $\Pi=\C \times 0^{n-1}$.  Since
the matrices defining $\Lambda_1$ and $\Lambda_2$
agree below the first row, we see that the
two sets $\Lambda_1(\Pi)$ and $\Lambda_2(\Pi)$
are identical.  But this translates into the
statement that the leaf of $\cal F$ through the
origin is preserved by the double lattice PET.
Finally, if one of the leaves is preserved, then
all the leaves are preserved.
\endproof

Let $\cal L$ denote the leaf of $\cal F$ through
the origin. The double lattice PET preserves
$\cal L$ and induces an action on this real
$2$-dimensional space.  Thus, even without knowing
that we have simply recreated the compactifications
from Theorem \ref{double2}, we can see that the
double lattice PETs here have associated planar
actions.  In the next section, we will
identify
these planar actions with the square turning maps
from Theorem \ref{double2}.

\subsection{Connection to Theorem \ref{double2}}
\label{dlp2}

We will work out the case $n=6$ explicitly.
This is a representative case.  The
cases $n=10,14,18,...$ follow the same
pattern.  We will treat the case $n=2$
separately, and from a different point of
view, in the next chapter.  In our
discussion,
we flip back and forth between $\C^6$ and $\R^{12}$ using
the identification
$(z_1,...,z_6) \leftrightarrow (x_1,y_1,...,x_6,y_6)$.
Our matrices are defined over $\C$, but $Q^{12}=[-1/2,1/2]^{12}$
is defined over $\R$.
Let
\begin{equation}
s_{ij}=\frac{s_i}{s_j}.
\end{equation}
Referring to the construction in \S 4, the matrix
$L_g$ turns out to be
\begin{equation}
\displaystyle
\left[\matrix{
-i & (\!-\!1\!+\!i) s_{21} & (1\!+\!i) s_{31} & (1\!-\!i) s_{41} & (\!-\!1\!-\!i) s_{51} & (\!-\!1\!+\!i) s_{61} \cr
(\!-\!1\!-\!i) s_{12} & i & (1\!+\!i) s_{32} & (1\!-\!i) s_{42} & (\!-\!1\!-\!i) s_{52} & (\!-\!1\!+\!i) s_{62} \cr
(\!-\!1\!-\!i) s_{13} & (\!-\!1\!+\!i) s_{23} & 2+i & (1\!-\!i) s_{43} & (\!-\!1\!-\!i) s_{53} & (\!-\!1\!+\!i) s_{63} \cr
(\!-\!1\!-\!i) s_{14} & (\!-\!1\!+\!i) s_{24} & (1\!+\!i) s_{34} & 2-i & (\!-\!1\!-\!i) s_{54} & (\!-\!1\!+\!i) s_{64} \cr
(\!-\!1\!-\!i) s_{15} & (\!-\!1\!+\!i) s_{25} & (1\!+\!i) s_{35} & (1\!-\!i) s_{45} &-i & (\!-\!1\!+\!i) s_{65} \cr
(\!-\!1\!-\!i) s_{16} & (\!-\!1\!+\!i) s_{26} & (1\!+\!i) s_{36} & (1\!-\!i) s_{46} & (\!-\!1\!-\!i) s_{56} & i
}\right].
\end{equation}
In general, the pattern is $4$-periodic, except for a suitable shift in
the indices of $s_{ij}$.

The matrices $I_6$ and $L_g$ represent the two lattices
$\Lambda_1$ and $\Lambda_2$.  Here $I_6$ is the
identity matrix.  We seek a matrix $M_1$ which
represents $X_1$ in the sense that $X_1=M_1(Q^{12})$.
Let $e_1,...,e_6$ be the standard basis vectors.
Looking at the proof of Lemma \ref{transversality}, we see that
the first row of $M_1^{-1}$ is $e_1$ and for
$k>1$ the
$k$th row of $M_1^{-1}$ is 
\begin{equation}
L_1^t \circ ... \circ L_{k-1}^t(e_k).
\end{equation}
Using this formula to compute $M_1^{-1}$, and then
taking inverses, we see that $M_1$ is the matrix
\begin{equation}
\left[\matrix{
1&0&0&0&0&0 \cr
(1-i) s_{12} & 1 & 0 &0 &0 &0 \cr
(1-i) s_{13} & (1-i) s_{23} & 1 & 0 &0 &0 \cr
(1-i) s_{14} & (1-i) s_{24} & (1-i) s_{34} & 1 &0 &0 \cr
(1-i) s_{15} & (1-i) s_{25} & (1-i) s_{35} & (1-i) s_{45} & 1 &0 \cr
(1-i) s_{16} & (1-i) s_{26} & (1-i) s_{36} & (1-i) s_{46} & (1-i) s_{56} & 1}
\right]
\end{equation}
Let $M_2=L_g M$.  The matrix data for our double lattice
PET is given by $(M_1,M_2,I_6,L_g)$.

Now we change coordinates.
Let $D_S$ and $\chi_1$ be as
in \S \ref{exp000}.
Let $A=\chi_1D_S$.  We compute

\begin{equation}
\label{up1}
AI_6=\chi_1D_S, \hskip 20pt
AL_g=\chi_2D_S, \hskip 20 pt
AM_1=L_1D_S \hskip 20 pt
AM_2=L_2D_S.
\end{equation}
Thus, the double lattice PET
$(A(X_1),A(X_2),A(\Lambda_1),A(\Lambda_2))$
is precisely the one constructed in
\S \ref{exp000}.
Note finally that
\begin{equation}
\label{dlp3}
I:=A L_g A^{-1} = J,
\end{equation}
where $J$ is the involution discussed in
\S \ref{fol000}.  This, $A$ carries
the invariant foliation for 
$(X_1,X_2,\Lambda_1,\Lambda_2)$ to the
invariant foliation discussed in \S \ref{fol000}.

\subsection{Extra Symmetry}

There is one additional symmetry we mention, though we
will not need this symmetry for any purpose.
The matrix
\begin{equation}
L_1'=\left[\matrix{
-1&+ 1&+1&+1&\dots \cr
-1&-i& 0&0&\dots\cr
0&+i&-1& 0&\dots\cr
0&0&+1&+i&\dots\cr
\dots&\dots&\dots&\dots&\dots}\right]
\end{equation}
is equivalent to $L_1$ and defines the
same lattices.  We compute that
\begin{equation}
\label{zoop}
\chi_1=\left[\matrix{
+i &0&0&0&0 \cr
0&-i&0&0&0 \cr
0&0&+1&0&0 \cr
0&0&0&+i&0 \cr
0&0&0&0&-1 & \cr
\dots&\dots&\dots&\dots&\dots}\right] L_1'.
\end{equation}
Let $A$ be the matrix listed.  We have
$A_{33}=iA_{22}$ and $A_{44}=iA_{33}$, etc.
The entry $A_{11}$ is special.  Thus,
if we exclude the top left entry of $A$,
the pattern along the diagonal has period $4$.

Let $[X_j]$ denote the lattice generated by the
sides of $X_j$. Algebraicaly, $[X_j]$ is just the
$\Z[i]$-span of the columns of $\chi_jD_S$.
Geometrically, Equation \ref{zoop} says that there is
an order $4$ isometry $J_{11}$ which carries 
$[X_1]$ to $\Lambda_1$.  More generally, there is
an order $4$ isometry $J_{ij}$ which carries
$[X_i]$ to $\Lambda_j$ for any pair of indices.
The case $n=2$, which we will treat specially
in the next chapter, has even more symmetry.
In this case, we can find a single order $4$
isometry $J$ which has
the action $[X_1] \to \Lambda_1 \to [X_2] \to \Lambda_2$.

We wonder if we can replace $L_1$ and
$L_2$ by different matrices so as to arrange
a similar situation in higher dimensions.
The most natural thing would be to let
$L_{\pm}^*$ denote the matrix obtained
from $\chi_{\pm}$ by multiplying the
top row by $i$.  In this case, we could
take $J(z_1,z_2,...,z_n)=(iz_1,z_2,...,z_n)$.
However, the quadruple $(\chi_1,\chi_2,L_1^*,L_2^*)$
turns out not to be a special system, and for
for random choices of $\{s_1,s_2,s_3,s_4,s_5,s_6\}$
we saw that the parallelotope $\chi_1D_S(Q^{12})$ is
not a fundamental domain for the lattice defined
by $L_1^*D_S$.  So, this attempt does not work.
We mention this because we think that our
construction is the simplest possible one which
will work.

\newpage

\section{The Octagonal PETs}
 
One should view this chapter as an elaboration
of the case $n=2$ from the previous chapter.
Here we will take a different point of view.
We are not sure if the cases $n=6,10,14,...$
can be treated in the same way we treat
the case $n=2$ here.  This chapter mostly
repeats material from [{\bf S1\/}].
In this chapter, we use paramaters
$\{1,s\}$ rather than $\{s_1,s_2\}$.  This is
our habit for the $2$-grid systems.

\subsection{The Reflection Lemma}

The {\it eigenlattice\/} of a parallelotope
is the lattice generated by the vectors 
parallel to the sides of the parallelotope.
Clearly a parallelotope is the fundamental
domain for its eigenlattice. 
A {\it reflection\/} in a face of the
parallelotope $P$ is an order $2$
linear isometry whose fixed set is
a subspace parallel to a face of $P$.
The face in question need not 
be a top-dimensional face.

\begin{lemma}
Let $P$ be a parallelotope and $\Lambda$ be its
eigenlattice.  Let $I$ be a reflection in a
face of $P$. Then $P$ is a fundamental domain
for $I(\Lambda)$.
\end{lemma}

\startproof
Let $Q \subset R$ be spaces and $L$ a lattice.
We call $Q$ an
{\it overdomain in\/} $R$ for $L$
if, for any $p \in R$, there is some $\lambda \in L$
such that $p+\lambda \in Q$.
Since the covolume of $I(\Lambda)$
equals the volume of $P$.  It suffices to prove that
$P$ is an overdomain in $\R^n$ for $I(\Lambda)$.

Let $\Pi_+$ be the fixed space of $I$ and
let $\Pi_-$ be the orthogonal space.  By
construction $\Pi_-$ is the $(-1)$ eigenspace of $I$.
Let $\pi_-$ be orthogonal projection onto $\Pi_-$.
Since multiplication by $(-1)$ is an automorphism
of $\Lambda$ and $\Pi_-$ is the $(-1)$ eigenspace,
we have $\pi_-(I(\Lambda))=\pi_-(\Lambda)$.
Clearly $\pi_-(P)$ is an overdomain 
in $\Pi_-$ for $\pi_-(\Lambda)$.
Hence $\pi_-(P)$ is an overdomain in
$\Pi_-$ for
$\pi_-(I(\Lambda))$.  Hence
$Q=\pi_{-1}(\pi_-(P)$
is an overdomain in $\R^n$ for $I(\Lambda)$.
So, we can find $\lambda_1 \in I(\Lambda)$ such
that $p+\lambda_1 \in P \cap \Pi_+'$, where
$\Pi_+'$ is some fiber of
$\pi_-$.

Since $\Pi_+$ is parallel to a face of $P$,
and $P$ is a parallelotope, 
$P \cap \Pi_+'$ is a translate of 
$P \cap \Pi_+$.  Call this the {\it translation property\/}.

Now, $P \cap \Pi_+$ is an overdomain in $\Pi_+$ for
$\Lambda \cap \Pi_+$.  But $I$ acts as the
identity on $\Pi_+$. Hence $P \cap \Pi_+$ is
an overdomain in $\Pi_+$ for $I(\Lambda) \cap \Pi_+$.
By the translation property,
$P \cap \Pi_+'$ is an overdomain in
$Q \cap \Pi_+'$ for for $I(\Lambda) \cap \Pi_+$.
Hence there is some
$\lambda_2 \in I(\Lambda)$ such that
$p+\lambda_1+\lambda_2 \in P$.
Hence $P$ is an overdomain in
$\R^n$ for $I(\Lambda)$
\endproof

\subsection{The Real Case}

Figure 6.1 shows a scheme for a $2$-dimensional
double lattice PET. 
$X_1$ and $X_2$ are the origin-centered
translates of $F_1$ and $F_2$
respectively. $\Lambda_1$ and $\Lambda_2$ are the
eigenlattices respectively of
$L_1$ and $L_2$.

\begin{center}
\resizebox{!}{2.5in}{\includegraphics{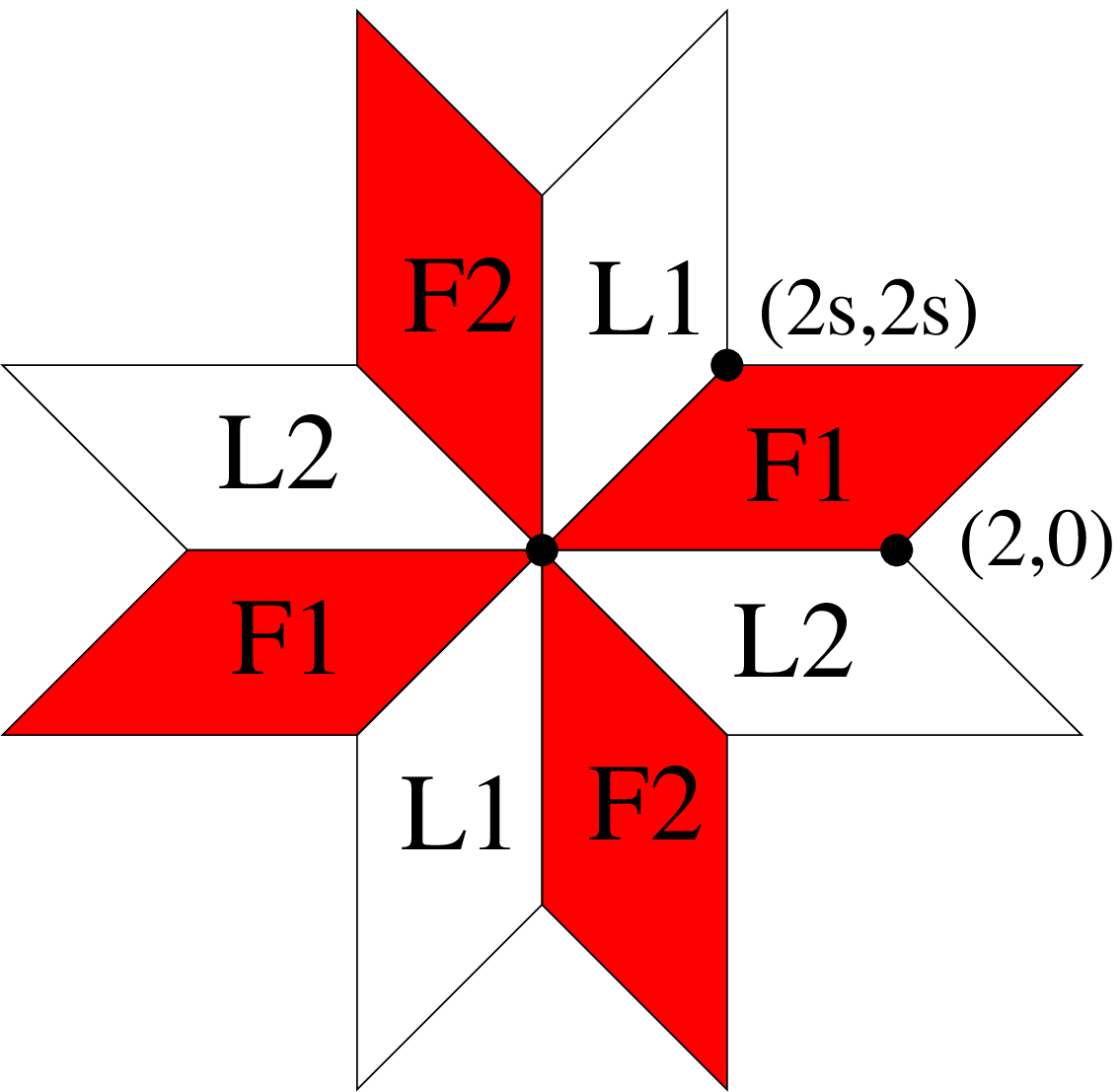}}
\newline
{\bf Figure 6.1:\/} Scheme for the octagonal PET.
\end{center}

Several applications of the Reflection Lemma
show that
$(X_1,X_2,\Lambda_1,\Lambda_2)$ 
really is the data for
a double lattice PET.
This double lattice PET depends on 
the parameter $s \in (0,1)$, which determines
the shapes of our parallelograms.
We call our system $(X_s,f_s)$.

In the next chapter we will give an account
of the following theorem, which we proved
(among many other things) in [{\bf S1\/}].

\begin{theorem}
\label{export}
Let $s \in (0,1)$ be irrational.  Then
almost every point of $(X_s,f_s)$ is periodic.
The sequence of periods is unbounded, and
the periodic islands are all semi-regular
octagons or squares. For almost every $s \in (0,1)$,
the system $(X_s,f_s)$ has uncountably
many aperiodic points.
\end{theorem}

\noindent
{\bf Remark:\/}
It would seem that our last statement of 
Theorem \ref{export} is not as strong as
it might be.  Since there are infinitely
many periodic islands and $X_s$ is compact,
there must be accumulation points of these
islands.  It would seem that these accumulation
points are aperiodic.  However, the issue
is that $f_s$ might not be defined at any
of these periodic points.  We were not
able to rule out this pathology for a
measure-zero set of parameters, though
we suspect aperiodic points exist for
all irrational parameters.

\subsection{The Complex Case}
\label{cc}

In this section, we
will define a $1$-parameter family of
$4$-dimensional lattice PETs 
$(X^{\C}_s,f^{\C}_s)$ which contains
$(X_s,f_s)$ as a $2$-dimensional invariant slice.
Our notation is meant to suggest the idea that
we produce these PETs by complexifying the
octagonal PETs.

Before we start, we remind the reader that
a {\it real plane\/} in $\C^2$ is a $2$-plane
$\Pi$ such that $\Pi$ and $i\Pi$ are orthogonal.
A {\it complex line\/} in $\C^2$ is a $2$-plane $\Pi$
such that $i\Pi$ and $\Pi$ are parallel
A {\it complex foliation\/}
is a $2$-dimensional foliation whose tangent
planes are complex lines.

We let $G$ denote the order $8$ dihedral
group generated by $\rho_1$ and $\rho_2$.
Here
\begin{equation}
\rho_1(z_1,z_2)=(-z_1,z_2), \hskip 30 pt
\rho_2(z_1,z_2)=(z_2,z_1).
\end{equation}
The invariant subspaces of $\rho_1$ and $\rho_2$
respectively are the complex lines given by
$\{z_2=0\}$ and
$\{z_1=z_2\}$. 
The group $G$ acts isometrically on $\C^2$. The
restriction of $G$ to $\R^2$ gives the group
of symmetries of Figure 6.1.

Let $X_1^{\C}$ be the parallelogram centered
at the origin and spanned by the vectors
\begin{equation}
(0,2), \hskip 20 pt
(0,2i), \hskip 20 pt
(s\zeta,s\zeta), \hskip 20 pt
(s \overline \zeta,s\overline \zeta);
\hskip 40 pt \zeta=1+i.
\end{equation}
The first two vectors lie in the complex line
fixed by $\rho_1$ and the second two vectors
lie in the complex line fixed by $\rho_2$.
We let $X_2^{\C}=\rho_1 \circ \rho_2(X_1^{\C})$.
We let $\Lambda_1^{\C}$ be the eigenlattice for
$\rho_1(X_1^{\C})$.  We could equally well
describe $\Lambda_1^{\C}$ as the eigenlattice for
$\rho_2(X_2^{\C})$.
We let $\Lambda_2^{\C}$ be the eigenlattice for
$\rho_2(X_1^{\C})$.  We could equally well
describe $\Lambda_1^{\C}$ as the eigenlattice for
$\rho_1(X_2^{\C})$.   Note that
$\Lambda_k^{\C} \cap \R^2=\Lambda_k$.
Several applications of the Reflection Lemma
show that $(\X_1^{\C},X_2^{\C},\Lambda_1^{\C},
\Lambda_2^{\C})$ is a double lattice PET.

By construction, we have
\begin{equation}
X_k^{\C} \cap \R^2 = X_k, \hskip 30 pt
\Lambda_k^{\C} \cap \R^2 = \Lambda_k,
\hskip 30 pt k=1,2.
\end{equation}
Here $(X_1,X_2,\Lambda_1,\Lambda_2)$ are
the data for the octagonal PET at the
parameter $s$.  
Hence, the octagonal PET $(X_s,f_s)$ is an invariant
``real slice'' of $(X_s^{\C},f_s^{\C})$.

The complex octagonal PETs have quite a bit of symmetry.
Let $\Gamma$ be the order $8$ dihedral group generated by
the elements
\begin{equation}
\iota_1(z_1,z_2)=(\overline z_1,\overline z_2), \hskip 30 pt
\iota_2(z_1,z_2)=(iz_1,iz_2).
\end{equation}
Each element of $\Gamma$ preserves each
parallelotope and each lattice defined
in connection with the complex octagonal PETs.
Hence, $\Gamma$ acts as an order $8$
group of symmetries of a complex octagonal PET.

There are $4$ elements of $\Gamma$ which act as
real reflections -- i.e., they pointwise fix real planes
in $\C^2$.  The planes $\Pi_0=\R^2$ and $\Pi_2=i\R^2$
are two of the fixed planes.   The other two fixed
planes are 
$\Pi_1=\zeta \R^2$ and
$\Pi_3=\overline \zeta \R^2$.
More simply,
\begin{equation}
\Pi_k=\zeta^k \R^2, \hskip 30 pt k=0,1,2,3.
\end{equation}
For instance, the map $\iota_2 \circ \iota_1$ fixes
$\Pi_1$ pointwise.  
It turns out that
\begin{itemize}
\item 
$(X_s^{\C},f_s^{\C}) \cap \Pi_0$ is a copy
of the octagonal PET $(X_s,f_s)$. 

\item 
$(X_s^{\C},f_s^{\C}) \cap \Pi_2$ is copy
of the octagonal PET $(X_s,f_s)$.

\item 
$(X_s^{\C},f_s^{\C}) \cap \Pi_1$ is a copy
of the octagonal PET $(X_{s/2},f_{s/2})$. 

\item 
$(X_s^{\C},f_s^{\C}) \cap \Pi_3$ is a  copy
of the octagonal PET $(X_{s/2},f_{s/2})$. 
\end{itemize}

Here, the word {\it copy\/} means {\it up to
a similarity\/}.
We already derived the first of these assertions above,
and the second one is not hard to see from symmetry.
We sketch a proof of the last two assertions
in [{\bf S1\/}, \S 5].  We will not use them here,
but we point them out in order to highlight some
of the beautiful symmetry of the complex
octagonal PETs.  Each complex octagonal PET
contains copies of the real octagonal PETs at
two different parameter values!

At this point, the reader might wonder whether
the complex octagonal PET is somehow just a product
of two of the octagonal PETs it contains, the
slice in $\Pi_0$ and the perpendicular slice
in $\Pi_2$.  This is not the case.  When
$s=1$, the two lattices $\Lambda_1$ and
$\Lambda_2$ coincide, and both are equal
to the beautiful $E_4$ lattice.  On
the other hand, the product lattice when
$s=1$ would just be $\Z^4$.

In the last section of this chapter,
we will discuss more connections
between the complex octagonal PETs and
the $E_4$-lattice.  All in all, one might say
that the complex octagonal PETs relate to
the real octagonal PETs sort of in the way
that the $E_4$ lattice relates to $\Z^2$.

\subsection{Connection to Square Turning}

Theorem \ref{double2} produces a double lattice
PET for the word $g=(A_1A_2)^2$ and for any
irrational parameter $s$.  In [{\bf S1\/}] we
recognized this double lattice PET
as the complex octagonal PET at parameter $s$.
In this section we will repeat the arguments,
through a little more tersely.  

The case $n=2$ in \S \ref{dlp2} gives rise to
the matrix data $(M_1,M_2,I,L_g)$ for
the double lattice PET. 
\begin{eqnarray}
\nonumber
I=\left[\matrix{1&0\cr0&1}\right]_{\C}=\left[\matrix{
1&0&0&0\cr
0&1&0&0\cr
0&0&1&0\cr
0&0&0&1}\right]_{\R}\\
\nonumber
L_g=\left[\matrix{-i&-s+is\cr -1/s-i/s&i}\right]_{\C} 
=\left[\matrix{
0&-1&-s&s \cr 
1&0&-s&-s \cr 
-1/s & -1/s & 0&1 \cr
1/s & -1/s &-1 &0}\right]_{\R} \\
\nonumber
M_1=\left[\matrix{1&0\cr 1/s-i/s&1}\right]_{\C}=
\left[\matrix{
1&0&0&0 \cr
0&1&0&0 \cr
1/s&-1/s &1&0 \cr
1/s&1/s&0&1}\right]_{\R} \\
M_2=L_gM_1=\left[\matrix{i&-s+is \cr 0 & i}\right]_{\C}=
\left[\matrix{
0&1&-s&s \cr
-1&0&-s&-s \cr
0&0&0&1 \cr
0&0&-1&0}\right]_{\R}
\end{eqnarray}
Here we are writing each matrix in two ways,
as a matrix over $\C$, and as a matrix over $\R$.
We are using the identification
\begin{equation}
(x_1+iy_1,x_2+iy_2) \leftrightarrow (x_1,y_1,x_2,y_2).
\end{equation}

Now we transform the picture by a suitable real
linear transformation.
We introduce the matrix
\begin{equation}
\gamma = \left[\matrix{
-2&0&+s&+s \cr
0&0&+s&+s \cr
0&+2&-s&+s \cr
0&0&-s&+s}\right]
\end{equation}

We compute
\begin{eqnarray}
\label{marked1}
\nonumber
\gamma M_1=
\left[\matrix{
0&0&+s&+s\cr
+2&0&+s&+s\cr
0&0&-s&+s\cr
0&+2&-s&-s}\right] \\
\nonumber
\gamma M_2=
\left[\matrix{
0&-2&+s&-s\cr
0&0&+s&+s\cr
+2&0&-s&+s\cr
0&0&-s&+s}\right] \\
\nonumber
\gamma I=
\left[\matrix{
-2&0&+s&+s\cr
0&0&+s&+s\cr
0&+2&-s&+s\cr
0&0&-s&+s}\right] \\
\nonumber
\gamma L_g=
\left[\matrix{
0&0&+s&-s\cr
0&-2&-s&+s\cr
0&0&+s&+s\cr
+2&0&-s&-s}\right] \\
\end{eqnarray}
\begin{equation}
\label{INVO}
\gamma L_g \gamma^{-1}=
\left[\matrix{
0&0&0&-1\cr
0&0&1&0 \cr
0&1&0&0 \cr
-1&0&0&0}\right]
\end{equation}

Let $(X_1',X_2',\Lambda_1',\Lambda_2')$ be
the double lattice PET determined by the
above matrices.  To recognize this as
the complex octagonal PET, we
identify $\R^4$ with $\C^2$ in a different
way:
\begin{equation}
\label{ident}
(x_1,x_2,y_1,y_2) \to (x_1+iy_1,x_2+iy_2).
\end{equation}

\begin{lemma}
$X_1'=X_1^{\C}$.
\end{lemma}

\startproof
With this identification, the sides of
$X_1'=\gamma M_1(Q^4)$ are
\begin{equation}
\label{X1prime}
(0,2), \hskip 20 pt (0,2i), \hskip 20 pt
(s\zeta,s\zeta), \hskip 20 pt
(s\overline \zeta,s\overline \zeta).
\end{equation}
This is exactly $X_1^{\C}$.
\endproof

\begin{lemma}
\label{arg1}
$X_2'=X_2^{\C}$.
\end{lemma}

\startproof
The second paralellotope $X_2'=\gamma M_2(Q^4)$ is
the image of $X_1'$ under the map
$\gamma L_g \gamma^{-1}$.   Looking at Equation
\ref{INVO}, we see that with our new identification
of $\C^2$ and $\R^4$, the map
$\gamma L_g \gamma^{-1}$ has the same action as
the map $-i \rho_1 \circ \rho_2$.  But then
\begin{equation}
X_2'=-i \rho_2(\rho_1(X_1'))=
-i\rho_1\rho_2(X_1^{\C})=-iX_2^{\C}=X_2^{\C}.
\end{equation}
The last equality comes from the fact that
multiplication by $(-i)$
is a symmetry of $X_2^{\C}$.
\endproof

\begin{lemma}
$\Lambda_1'=\Lambda_1^{\C}$
\end{lemma}

\startproof
The lattice $\Lambda_1'$ represented by $\gamma I_4$ is the
$\Z$ span of the vectors
\begin{equation}
\label{L1prime}
(2,0), \hskip 20 pt
(2i,0), \hskip 20 pt
(s\overline \zeta,s\overline \zeta),
\hskip 20 pt
(s\zeta,s\zeta).
\end{equation}
technically, the first vector we read off is $(-2,0)$, but
changing the sign has no effect on the lattice.
Recall that $\rho_2(z_1,z_2)=(z_2,z_1)$.  If
we apply $\rho_2$ to the vectors listed in
Equation \ref{X1prime}, we get the vectors in
Equation \ref{L1prime}. But then
$\Lambda_1'$ is the image under $\rho_2$ of the
eigenlattice for $X_1^{\C}$.   Since
$\Lambda_1^{\C}$ has the same description,
we see that $\Lambda_1'=\Lambda_1^{\C}$.
\endproof

\begin{lemma}
$\Lambda_2'=\Lambda_2^{\C}$.
\end{lemma}

\startproof
This has the same proof as
Lemma \ref{arg1}.
\endproof

Thus we see that the double lattice PETs
produced by Theorem \ref{double2} in the
lowest dimensional case are the complex
octagonal PETs.  In other words,
the compactifications associated
to the alternating grid systems are the
complex octagonal PETs.  In the next section
we will deduce some dynamical consequences
from this fact, and from our knowledge
of the complex octagonal PETs.

\subsection{Dynamical Consequences}

Combining Theorem \ref{double2} with the
explicit calculation in the previous
section, we get the following corollary.

\begin{corollary}
\label{goodslice}
The complex octagonal PET
$(X_s^{\C},f_s^{\C})$ has two invariant
irrational orthogonal complex foliations,
$F_1$ and $F_2$.  Here
\begin{enumerate}
\item $F_1$ is spanned by 
$(\zeta,\overline \zeta)$ and
$(\overline \zeta,-\zeta)$.
\item $F_2$
is spanned by $(\overline \zeta,\zeta)$ and
$(\zeta,-\overline \zeta)$.
\end{enumerate}
For any $z \in \C$, there is a leaf $L=L(s,z)$ of
$F_1$ (or of $F_2)$ such that the
restriction of $f_s^{\C}$ to $L$
is conjugate
to the action of $g_{s,z}$ on $\C$.
The map $C \to L$ is a piecewise
isometry relative to the Euclidean
metric on $\C$ and the path metric on $L$.
\end{corollary}

\noindent
{\bf Remarks:\/}
\newline
(i)
We get a piecewise isometry from
$\C$ to $L$ in our corollary because
$X^{\C}$ is a parallelopiped rather than a torus.
Were we to glue together the opposite sides of
$X^{\C}$, the map $C \to L$ would be an
isometry.
\newline
(ii) 
One surprising thing about our corollary is that
the complex octagonal PETs have two invariant
foliations, and the restriction of the map to
a leaf in either foliation is conjugate to the
square turning map associated to
$(A_1A_2)^2$.  
\newline

\noindent
{\bf Proof of Theorem \ref{stationary1}:\/}
Let $s$ be some irrational parameter and
let $z \in \C$ be arbitrary.  Let $L$ be
as in Corollary \ref{goodslice}.
Let $\{p_k\}$ be the unbounded sequence of periods
associated to the real octagonal PET $(X_s,f_s)$, as
guaranteed by Theorem \ref{export}.
For any $k$, there is some periodic island in
$(X_s,f_s)$ of period $p_k$.   But
$(X_s,f_s)$ is an invariant slice of
$(X_s^{\C},f_s^{\C})$ and this latter system is a
PET.  Hence, the points of our $2$-dimensional
periodic island are contained in a $4$-dimensional
periodic island $P$ of period $p_k$. 
 
The piecewise isometry from $\C$ to 
$L$ induces a locally Euclidean measure to $L$,
relative to which
$L \cap P$ has positive
density in $L$.  This density is exactly
${\rm vol\/}(P)/{\rm vol\/}(X)$.  Hence,
a positive density set of points in the
invariant leaf $L$ have period $p_k$.
But then the same statement holds for points
in $\C$, relative to the map $g_{z,s}$.
Since $g_{s,z}$ has a positive density set
of points of period $p_k$, it has infinitely
many $p_k$-periodic islands.
\endproof

\noindent
{\bf Proof of Theorem \ref{stationary2}:\/}
We keep the same notation as in the proof
of Theorem \ref{stationary1}.  This time,
we choose $s$ so that the octagonal PET
$(X_s,f_s)$ has uncountable many 
aperiodic points.

A routine calculation shows that the plane spanned by
 $(\overline \zeta,\zeta)$ and
$(\zeta,-\overline \zeta)$ is transverse to $\R^2$.
Hence, the leaf $L$ intersects $\R^2$ in a
countable collection of points.  The same
goes for any leaf of $F_1$. Hence
there uncountably many leaves of $F_1$
which contain aperiodic points of $X \subset \R^2$, the
domain for the real octagonal PET.  As in 
Lemma \ref{count}, this means that there are
uncountably many choices of $z \in \C$ such
that $L=L(s,z)$ contains an aperiodic point
of the real octagonal PET.

Let $\omega: \C \to L$ be the conjugacy
guaranteed by Corollary \ref{goodslice}.
Let $z \in \C$ be such that $L$ contains an aperiodic
point $x \in L \cap X$.
The orbit of $O(x)$ of $x$ is infinite, contained
in $X$, and also contained in $L$.
Let $x'=\omega^{-1}(x)$.  By
construction the orbit $O(x')$ under
$g_{s,z}$ is infinite.  Moreover
$\omega$ carries $O(x')$ to $O(x)$.
We want to see that $O(x')$ is unbounded.

If $B$ is any bounded subset of $\C$, then
$\omega(B)$ intersects $X$ in only
finitely many points.  This follows from
the fact that $\omega$ is an isometry when 
$X^{\C}$ is interpreted as a torus, and
the image $L=\omega(\C)$ is transverse to $X$.
Since $\omega(B)$ can only contain finitely
many points of $O(x)$, the bounded
set $B$ cannot contain all of $O(x')$.
Hence $O(x')$ is unbounded.
\endproof

\newpage

\section{Renormalization and its Self-Similarity}

\subsection{Renormalization}

Theorem \ref{export} derives from a 
{\it renormalization scheme\/} we found for
the (real) octagonal PETs.
Let $(X_s,f_s)$ be the octagonal PET at
parameter $s$.  
Given $Y \subset X$ (and suppressing the parameter) let
$f|Y$ denote the first return map of $f$ to $Y$,
assuming that this map is well-defined.

Define the map
$R: (0,1) \to [0,1)$ by the formula
\begin{itemize}
\item $R(s)=1-s$ if $s \in [1/2,1)$.
\item $R(s)=1/(2s)-{\rm floor\/}(1/(2s))$ if
$s \in (0,1/2)$.
\end{itemize}
For all but countably many choices of $s$, we have
$t=R(s)>0$.  In all these cases, we prove
the following result.
\begin{theorem}
\label{renorm}
There are clean convex polygons
$Y_t \subset X_t$ and $Z_s \subset X_s$ and
a similarity $\phi_s: Y_t \to Z_s$ such that 
\begin{itemize}
\item $\phi_s$ conjugates $f_t|Y_t$ to $f_s^{-1}|Z_s$.
\item Every nontrivial orbit of $f_t$ intersects $Y_t$ and
every nontrivial orbit of $f_s$, escept possibly
for certain orbits of period $2$, intersects $Z_s$.
\end{itemize}
The map $\phi_s$ is an isometry when $s>1/2$ and
a contraction when $s<1/2$.
\end{theorem}
A polygon is {\it clean\/} if its boundary does
not intersect any of the open periodic islands.
We describe $Y$ and $Z$ explicitly in [{\bf S1\/}].

For any given parameters $s$ and $t=R(s)$, the
proof of Theorem \ref{renorm} is a fairly easy
and finite calculation.  Indeed, the way
we will establish Theorem \ref{export} for
the single parameter $\sqrt 2/2$ is just
to appeal to a finite calculation like this.
See below for details.

Establishing Theorem \ref{renorm}
for all parameters simultaneously is
also a finite calculation, but it is much
more involved. The idea is that we consider
a $3$-dimensional piecewise affine system
whose $2$-dimensional fibers are the octagonal
PETs. We then prove by direct calculation a
version of Theorem \ref{renorm} for the
$3$-dimensional system and observe that the
$3$-dimensional result reduces to the
$2$-dimensional result above in each slice.
All this is easier said than done, however.

\subsection{Sketch of Theorem \ref{export}}

When $s$ is irrational, the infinite sequence
$\{R^n(s)\}$ exists.  It follows almost
immediately from Theorem \ref{renorm} that
the octagonal PET $(X_s,f_s)$ has infinitely
many periodic islands when $s$ is irrational.
But this system can only have finitely many
periodic islands less than any given period.
Hence, the sequence of periods of points
in $(X_s,f_s)$ is unbounded.
This gives the first statement of Theorem
\ref{export}.

For the second statement, say that 
$p \in X_s$ is a {\it limit point\/} if
every open neighborhood of $p$ contains
infinitely many periodic islands.  Say
that $p$ is {\it bad\/} if the orbit of
$p$ is undefined.  Otherwise say that
$p$ is {\it good\/}.  Since $X_s$ is compact and
there are infinitely many periodic islands,
we know that there is at least one limit
point.  The tricky part is showing that
there are some good limit points.

We prove the following two statements in [{\bf S1\/}].
\begin{enumerate}
\item The set of
limit points of $(X_s,f_s)$ has positive
$1$-dimensonal Hausdorff measure for
all irrational $s$.
\item The set of bad limit
points has zero $1$-dimensional Hausdorff
measure for almost all $s$.
\end{enumerate}
This leaves
some good limit points for almost all $s$.
We will sketch the proof of the Statement 1,
 because we especially like the
argument, and because it is related to the
discussion in \S \ref{E4symm}.

Let $D_4$ be the order $8$ dihedral symmetry
group of the unit square $Q^2$.  We say that a
dihedral {\it polygon\/} is a polygon which, up
to translation, has $D_4$ as its symmetry group.
To be clear, if we were to rotate $Q^2$ by
a typical angle, it would not be a dihedral polygon.
The dihedral polygons are either squares or
semi-regular octagons.

We check, for the parameters $s=1/(2n)$, that all
the periodic islands are dihedral polygons and
right-angled isosceles triangles.  We then use
Theorem \ref{renorm} and induction to deduce
the same result for all rational parameters $s$.
Moreover, we use the explicit scaling factors
in Theorem \ref{renorm} to show that
the diameter of any triangular
tile is at most $2^{2-k/4}$, where $k$ is
the length of the $R$-orbit of $s$. This
quantity decays exponentially with $k$. Taking
a limit, we see that every periodic island 
in the irrational case is a dihedral polygon.
 Figure 7.1 shows a picture of the tiling
for the rational parameter $s=13/21$.  The
yellow triangular tiles are quite small.

\begin{center}
\resizebox{!}{2.4 in}{\includegraphics{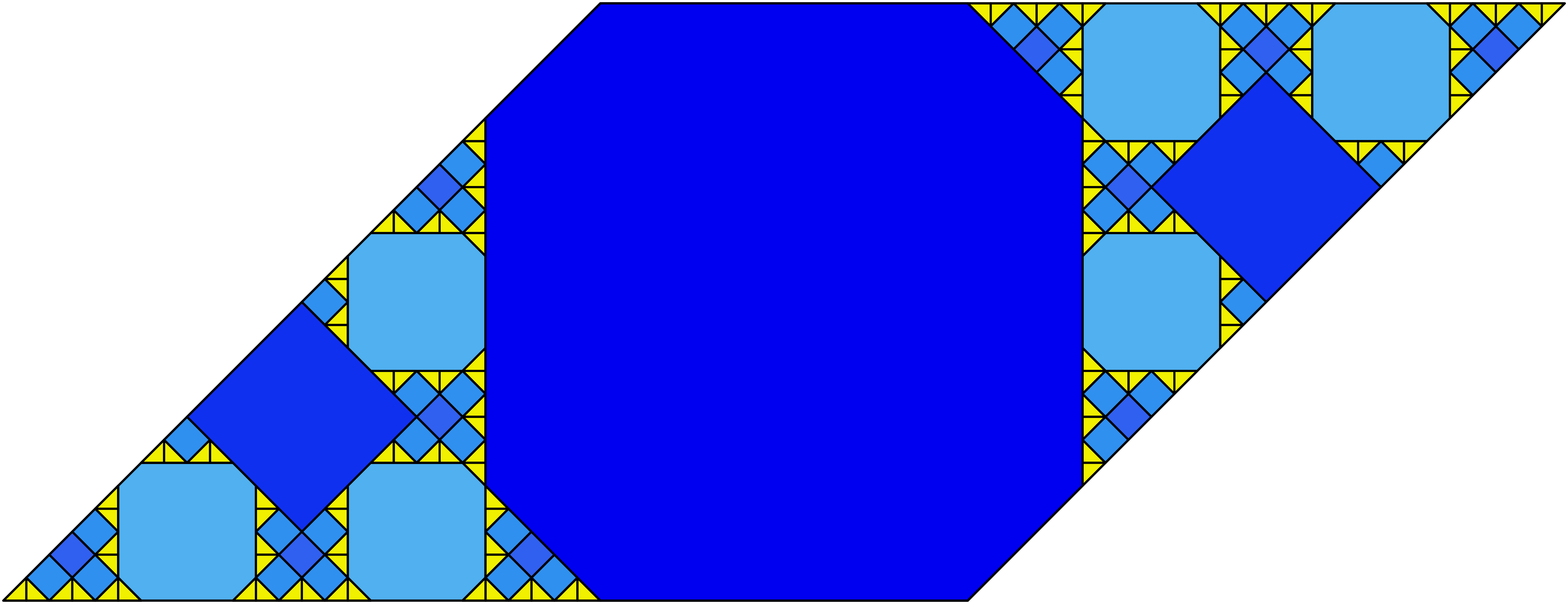}}
\newline
{\bf Figure 7.1:\/} The tiling for $s=13/21$.
\end{center}

There are $3$ kinds of edges in a $D_4$ polygon:
horizontal, vertical, and diagonal.  Consider how
a horizontal line $L$ intersects a dihedral polygon.
Assuming that $L$ does not contain a vertex of $P$,
then $L$ either intersects $P$ in two diagonal
sides or in two vertical sides.  Moreover, if
two $D_4$ polygons are tangent along an edge,
and $L$ intersects the interior of this edge,
then $L$ intersects both polygons in the
same kinds of edges, either horizontal or
vertical.

In all cases, the island of period $1$ -- i.e.,
the {\it fixed island\/} -- is just $X_1 \cap X_2$.
When $s<1/2$ this set is a square with sides
parallel to the coordinate axes.  In Figure
7.1, this is the half-shown big red square
on the right.  Consider the case $s<1/2$ for
ease of exposition. 
All the horizontal lines of the form $y=y_0 \in [-s,s]$
intersect $X_s$.  If we throw out countably
many choices of $y_0$, then the remaining lines
do not contain vertices of periodic islands.

The horizontal line  $L$ intersects the left
edge of $X_s$ in a diagonal edge.  On the other
hand, $L$ intersects the fixed island in a vertical
edge.   Given what we have said about how $L$
interacts with semi-regular octagons, we
see that $L$ cannot simply run through the
interior of a finite union of periodic islands.
That is, $L$ must contain a limit point.
This means that the projection of the set of
limit points onto the vertical axis contains
an interval.  Hence, the set of limit points
has positive $1$-dimensional Hausdorff measure.

The projection result is rather surprising.
The set of limit points in Figure 7.2 below
is a Cantor set. Nonetheless, its vertical 
projection contains an interval.

\subsection{Self-Similar Examples}

Figures 7.2 and 7.3 shows the picture respectively
for the two related parameters
$s=\sqrt 2/2$ and $s=\sqrt 2/4$.
Both pictures are self-similar and both figures
arise as different slices of the complex
octagonal PET $(\X_s^{\C},f_s^{\C})$ at the
parameter $s=\sqrt 2/2$.

\begin{center}
\resizebox{!}{2in}{\includegraphics{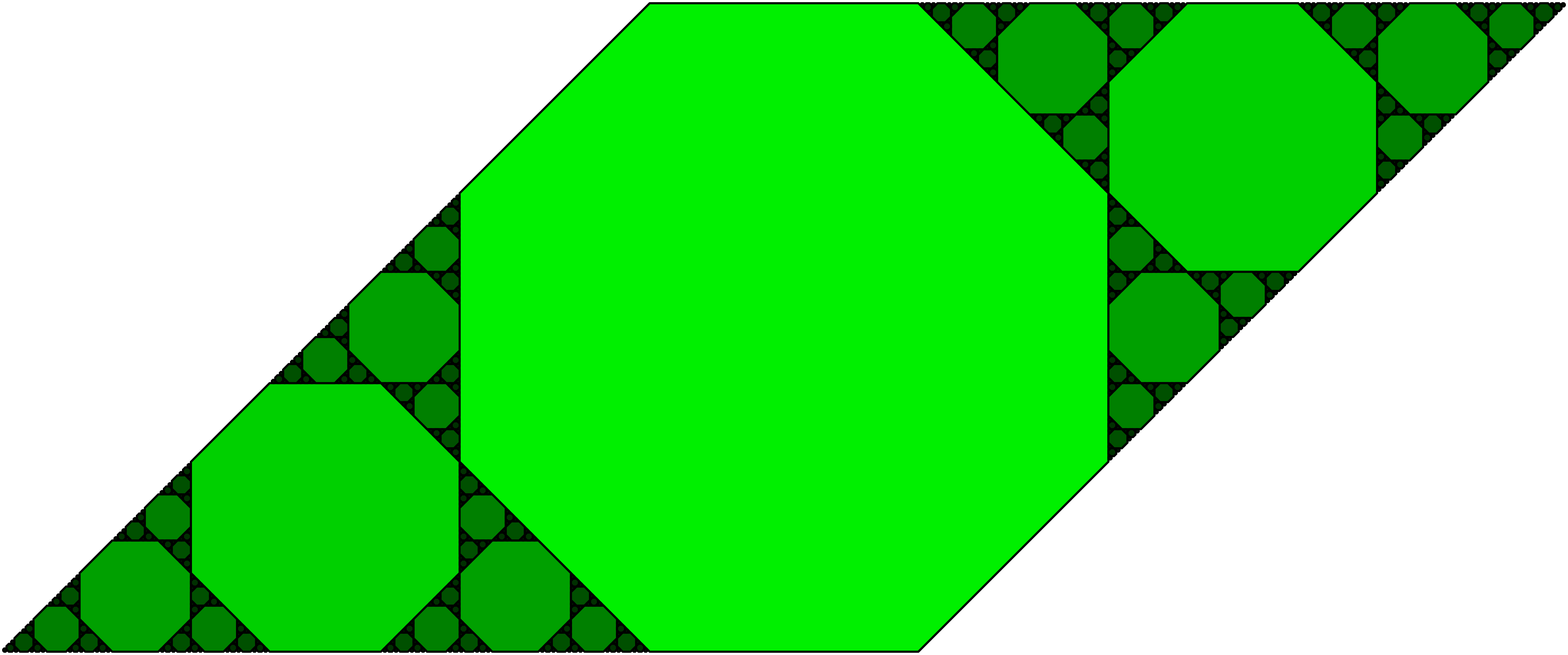}}
\newline
{\bf Figure 7.2:\/} The tiling for $s=\sqrt 2/2=0:1:2:2:2...$.
\end{center}

\begin{center}
\resizebox{!}{2.5 in}{\includegraphics{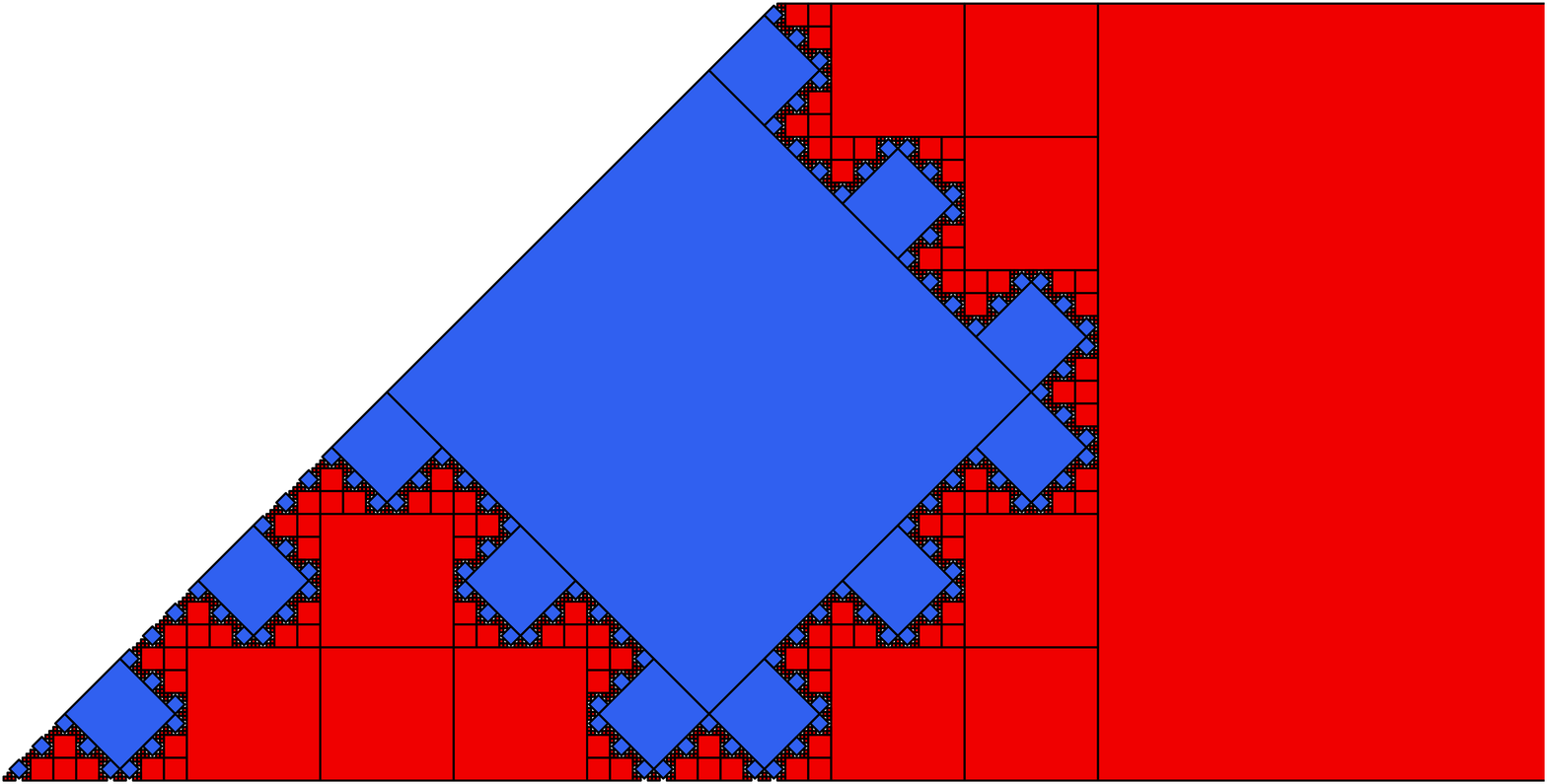}}
\newline
{\bf Figure 7.3:\/} The tiling for $s=\sqrt 2/4=0:2:1:4:1:4...$
\end{center}

These particular cases are very similar to
other systems which arise in this kind of dynamics.
  Figure 7.2 is locally isometric
to the tiling produced by the main example in
[{\bf AKT\/}], and also to the tiling produced
by outer billiards on the regular octagon.
Figure 7.3 is locally isometric to one of the
tilings produced by the Truchet tile system in
[{\bf Hoo\/}].

Figure 7.4 shows another example.
Another one of our result from
[{\bf S1\/}] is that the periodic tiling of
$(X_s,f_s)$ consists entirely of squares
if and only if the continued fraction expansion
of $s$ has the form $0:a_1:a_2:a_3...$ with
$a_k$ even for all odd $k$.  This condition
turns out to be equivalent to the condition
that $R^k(s)<1/2$ for all $k$.  We call
such parameters {\it oddly even\/}.

\begin{center}
\resizebox{!}{2.3in}{\includegraphics{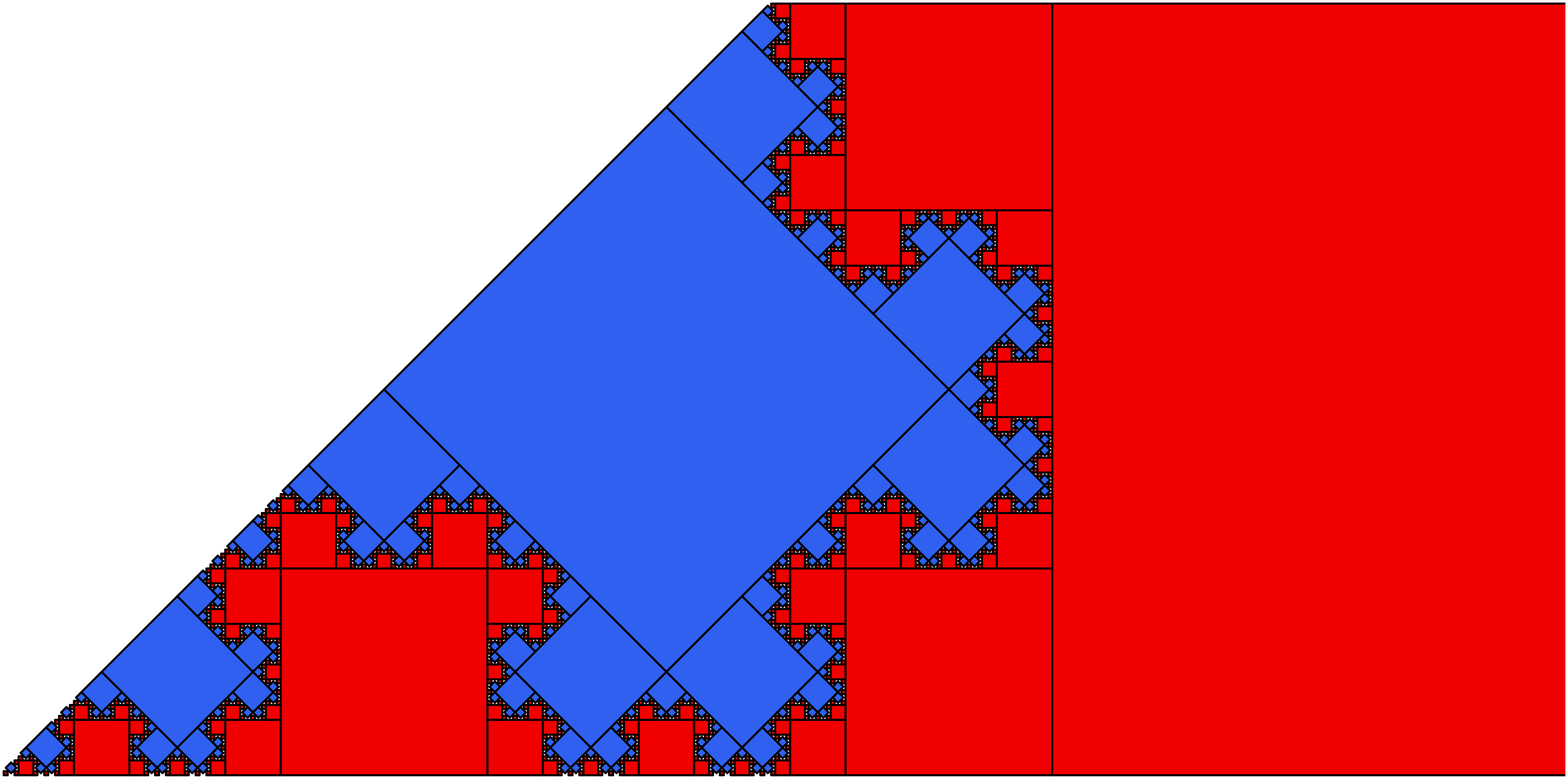}}
\newline
{\bf Figure 7.4:\/} The tiling for $s=\sqrt 3/2-1/2=:0:2:1:2:1:2...$.
\end{center}

In all three cases shown, the self-similar nature of
the picture derives from the fact that
$s$ is a periodic point of the
renormalization map $R$.
In these cases, a finite calculation shows that
the system $(X_s,f_s)$ is renormalizable, and
then one can deduce the structure of the tiling.
We will explain this in somewhat more detail for
the most familiar of the pictures, Figure 7.2.

In the cases shown here, 
the set of limit points is a fractal
having Hausdorff dimension greater than $1$.
For instance, the Cantor set in Figure 7.2,
corresponding to $s=\sqrt 2/2$, has dimension
$\log(3)/\log(1+\sqrt 2)$.  Since the bad set has
Hausdorff dimension at most $1$, there are
always good limit points.  

Now we describe something fairly amazing.
One could say that the existence of Figure 7.4 implies
the existence of some unbounded orbits for the square
turning system $g_{s,z}$ for $s=\sqrt 3/2-1/2$ and
a suitable choice of $z$.  Let $U$ be an unbounded
orbit corresponding to one of these limit points.
It seems that the rescaled limit of $U$ is locally
isometric to the fractal curve in Figure 7.4! 
Thus, the fractal in Figure 7.4 as
a kind of bird's eye view of a particular unbounded 
square-turning orbit.
The same statement seems to hold for
any oddly even parameter.  For
oddly even parameters, the square turning system seems to
``implement'' the fractal limit set of the
associated octagonal PET.  We have not yet tried for a proof.

\subsection{More Details in one Case}
\label{moredetail}

Here we discuss the case $s=\sqrt 2/2$ in more detail.
We set $f=f_s$, etc.

Figure 7.5 shows the partition of definition for
$f$ and Figure 7.6 shows the partition of definition
for $f^{-1}$.  Each polygon in Figure 7.5 is
translation equivalent to a unique polygon in Figure 7.6.
The map $f$ simply performs the translations which carry
each polygon in Figure 7.5 to the corresponding polygon
in Figure 7.6.  Thus, one can see the action of $f$
just by comparing the two figures.

\begin{center}
\resizebox{!}{2.3in}{\includegraphics{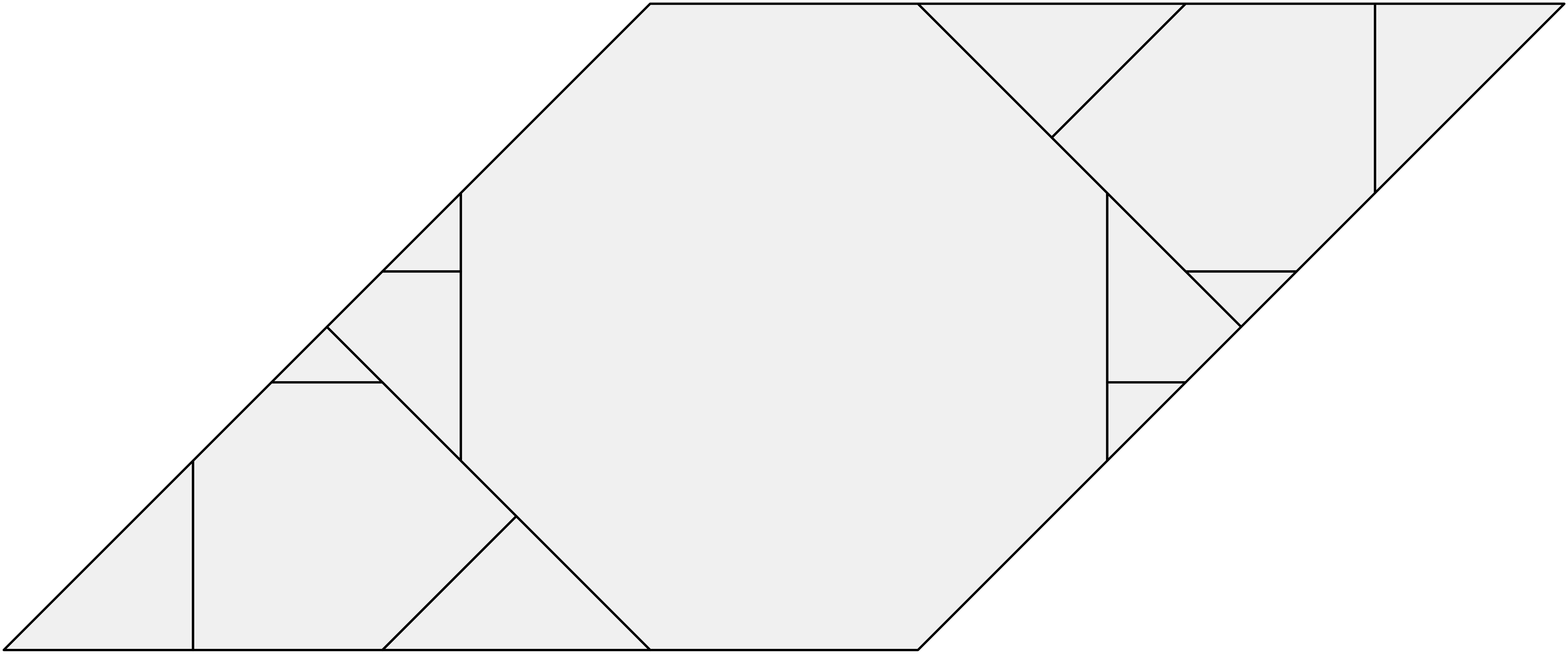}}
\newline
{\bf Figure 7.5:\/} The forward partition for $s=\sqrt 2/2$.
\end{center}

\begin{center}
\resizebox{!}{2.1in}{\includegraphics{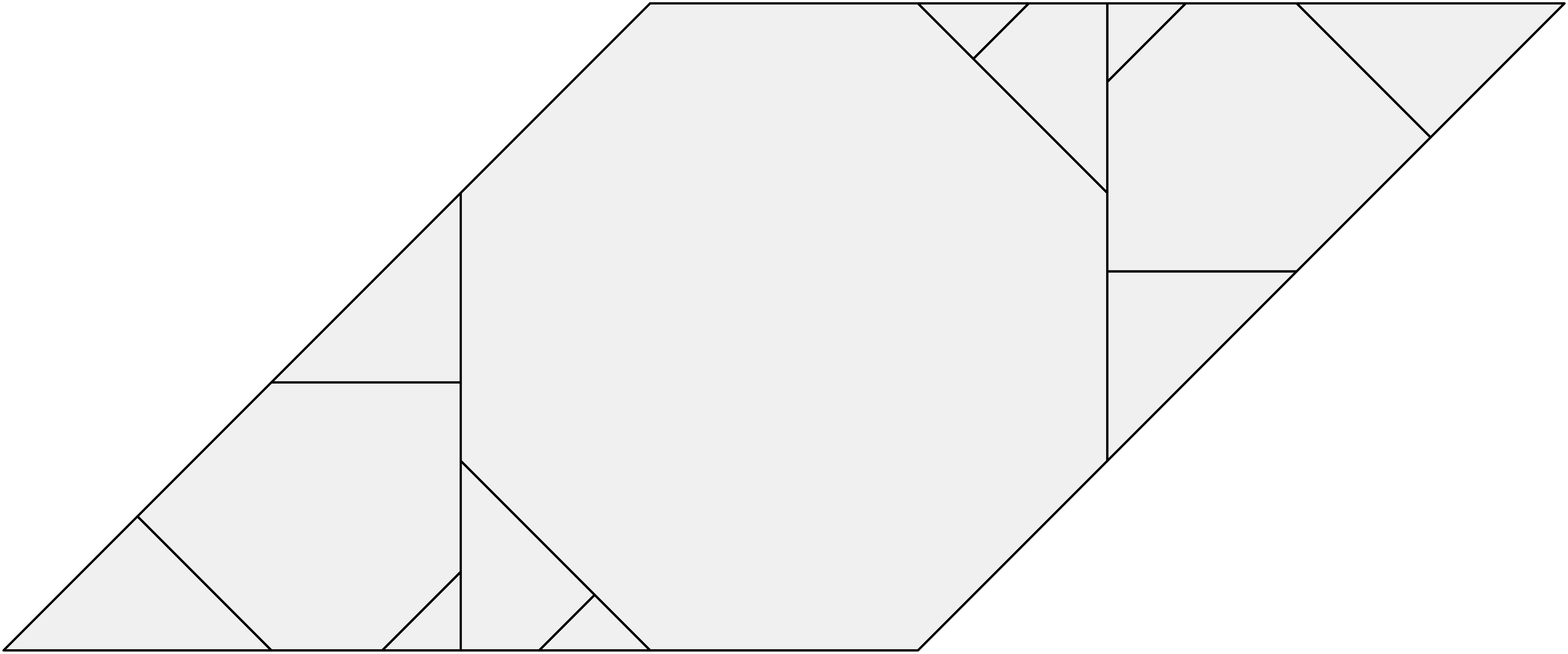}}
\newline
{\bf Figure 7.6:\/} The backward partition for $s=\sqrt 2/2$.
\end{center}

Inspecting the figures, we see that the central octagon $O_1$ is 
the fixed island of $f$ and that there are two smaller
octagonal islands $O_{21}$ and $O_{22}$
of period $2$, as shown in Figure 7.7.

\begin{center}
\resizebox{!}{2.3in}{\includegraphics{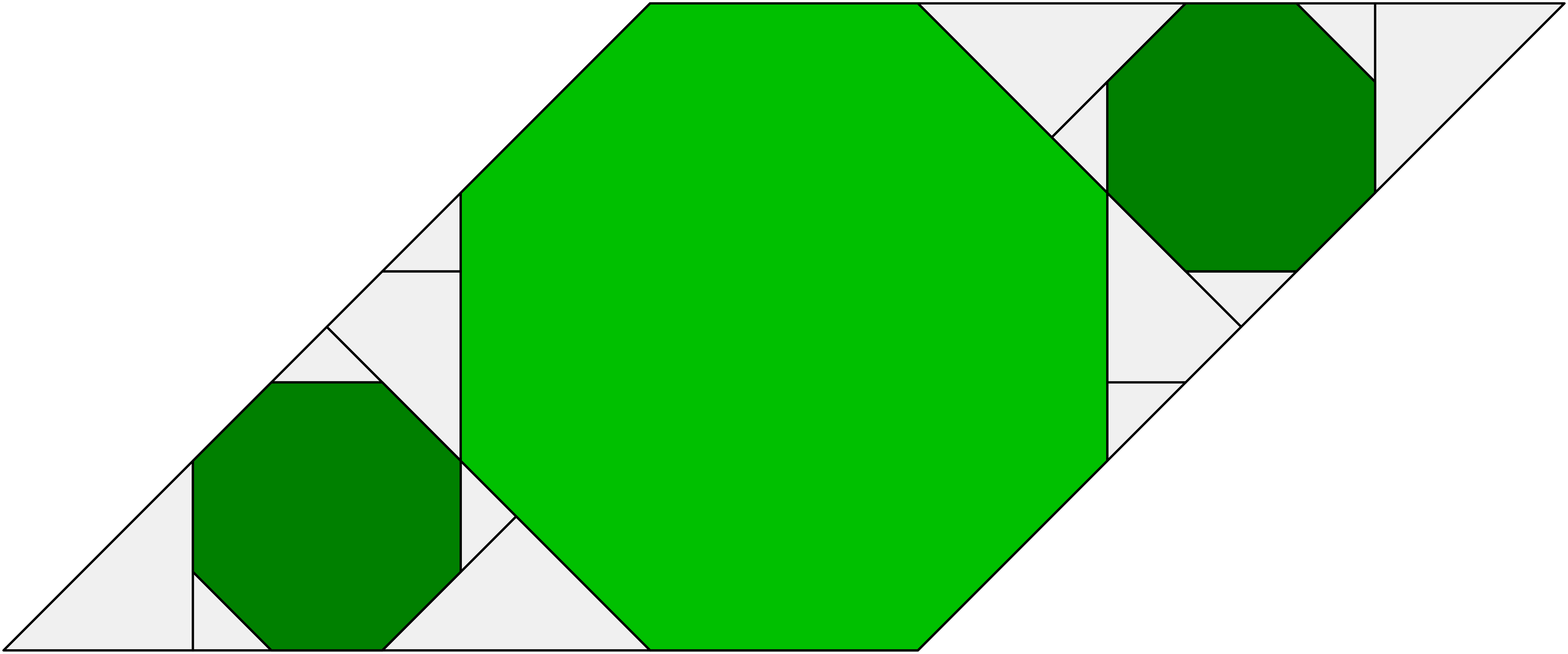}}
\newline
{\bf Figure 7.7:\/} The period islands $O_1$ and $O_{21}$ and $O_{22}$.
\end{center}

There are $8$ sets of interest to us.  $X-O_1$ is a union
of two kites $K_{11}$ and $K_{12}$ and
$X-O_1-O_2$ is a union of $6$ smaller kites
$K_{21},..,K_{26}$.  There are two similarities carrying
$K_{1j}$ to $K_{2k}$ for $j \in \{1,2\}$ and $k \in \{1,...,6\}$.
One can check with a finite calculation that each
of these $12$ maps conjugates $f|K_{1j}$ either to $f|K_{2k}$
or $f^{-1}|K_{2k}$.  This property implies that the tiling
by periodic islands is invariant under the action of these
$12$ maps.  The rest of Figure 7.2 is then determined by
this structure.

We have reduced the proof of Theorem \ref{export},
for the parameter $s=\sqrt 2/2$, to $12$ calculations.
We can use symmetry to cut down on the amount of work
we have to do.
Let $K_{11}$ be the big kite on the right and
let $K_{21}$ be the rightmost small kite.
Let $K_{22}$ be the bottom kite on the right.
The left edge of $K_{22}$ is the right edge of $O_1$.

Let $\rho_1$ denote reflection through the origin.
Let $\rho_2$ denote reflection through the long
diagonal of $X$.  It follows from the definition of
the octagonal PETs that $f$ commutes with $\rho_1$.
Interchanging the roles played by the two sides of
$X_1$ in the definition of the system, we see that
$\rho_2$ conjugates $f$ to $f^{-1}$.
Using this symmetry, we see that it suffices to check
just $2$ of the $12$ conjugacies mentioned above, 
namely one of the maps $K_{11} \to K_{12}$ and
one of the maps $K_{11} \to K_{22}$.

We will discuss the calculation for $K_{11} \to K_{21}$
in the next section.  The calculation $K_{11} \to K_{22}$
can be done in the same way, though we
did not actually make the second calculation.  Rather, in
[{\bf S1\/}] we established a general symmetry
which reduces the second calculation to the first one.
See \S \ref{E4symm} for a discussion of this extra symmetry.

\subsection{The Calculation}

Define
\begin{equation}
K_j=K_{j1} \cup \rho_1(K_{j1}), \hskip 30 pt j=1,2.
\end{equation}
Let $\phi$ be the piecewise orientation-reversing
similarity which carries the left (respectively right)
half of $K_1$
to the left (respectively right)
half of $K_2$.   
If one can verify that
\begin{equation}
\label{sim}
\phi^{-1} \circ (f^{-1}|K_{2}) \circ \phi=f|K_{1}
\end{equation}
then the same relation holds with
$K_{j1}$ in place of $K_j$.

It turns out that there is a more precise
relationship in this case.
One can check that $K_2$ is
an invariant set for $f^{-3}$ and that
$\phi$ conjugates the action of $f$ on
$K_1$ to the action of $f^{-3}$ on $K_2$.
We close this chapter by sketching how
one makes this verification.

Let $\Pi_+$ denote the forward partition for $f$, shown
in Figure 7.5. 
Say that a finite sequence of points $p_1,...,p_k$ is
{\it feasible\/} if there are polygons
$P_1,...,P_k$ in the partition $\Pi_+$ such that
$p_j \in \overline P_j$ for all $j$ and $f_j(p_j)=p_{j+1}$
for $j=1,...,k-1$.  Here $f_j$ is the extension
of $f|P_j$ to the closure $\overline P_j$.
A finite portion of a genuine orbit is feasible,
and so are the limits of such things.
However, a feasible sequence may not be a portion
of a well-defined orbit.
We call another feasible sequence
$q_1,...,q_k$ {\it compatible\/} with
$p_1,...,p_k$ if $q_j \in \overline P_j$
for all $j$. That is, both sequences visit
the same sequence of partition polygons.
Here is a helpful lemma.

\begin{lemma}[Definedness Criterion]
\label{partition}
Suppose that $Q$ is a $n$-gon, with
vertices $p_{11},...,p_{n1}$.  If
there exist $n$ mutually compatible
sequences $p_{j1},...,p_{jk}$ for
$j=1,...,n$, then $f^k$ is well
defined on all points on the interior of $Q$.
\end{lemma}

\startproof
Let $P_1,...,P_k$ be the sequence of polygons
of $\Pi_+$ visited by our sequences.  Let
$q_1$ be a point in the interior of $Q$.
$q_1$ lies in the interior of $P_1$ by
convexity. Hence $f$ is defined on $q_1$
and $q_2=f(q_1) \in P_2$.  And so on.
\endproof

Let $g$ be the map on the left hand side of
Equation \ref{sim}.
Using the definedness criterion, with respect
to the inverse map $f^{-1}$ and the inverse
partition $\Pi_-$, we check that
$f^{-3}$ is well defined on every polygon of
the form $\phi(\Pi_+|K)$.
But this means that both $f$ and $g$ are
well defined (and hence translations)
on the interior of each
of the $12$ polygons of $\Pi_+-O_1$.
But then it suffices to check Equation
\ref{sim} on $12$ points, one per polygon.
We omit the details of these few calculations.

\newpage

\section{Connections to $E_4$}

\subsection{Three Connections}

There is a big literature on things related to $E_4$.
See, for instance [{\bf CS\/}].  To match what we
have done in previous chapters, we will scale
the $E_4$ lattice so that its shortest vectors have
length $2$.

The $E_4$ {\it lattice\/} $\Lambda$ is the lattice 
\begin{equation}
2\Z^4 \cup \Z^4_{\rm odd\/}.
\end{equation}
Here $\Z^4_{\rm odd\/}$ is the subset of
vectors having all odd coordinates.
Geometrically, $\Lambda$ is the union of
vertices and centers of the cubes in the
cubical grid of side length $2$ in $\R^4$.
What makes $\Lambda$ 
so symmetric is the geometric
miracle that the distance from the center of a
$4$-dimensional cube to a vertex of the cube
is the same as the side length of the cube.

The $E_4$-{\it polytope\/} $P$ is the convex hull
of the set of $24$  vectors of length $2$
in $E_4$.  Up to permuting the coordinates and/or
multiplying some of the coordinates by $-1$, these
vectors are all equivalent to $(2,0,0,0)$ or
to $(1,1,1,1)$.
The $E_4$ polytope
is one of the $4$-dimensional platonic solids.  It
enjoys $3$ properties.
\begin{enumerate}
\item $P$ is regular.  The symmetry group of $P$ acts
transitively on the complete flags of $P$.
\item $P$ tiles space.  The Voronoi cells of $\Lambda$
are all translates of the smaller copy $(1/2)P$.
\item $P$ is self dual.  Each length-$2$ normal to
a facet (i.e. codimension one face) of $P$ has the form
$(\sqrt 2,\sqrt 2,0,0)$ up to signs and
permutation.  The convex hull $P^*$ of these
unit normals is isometric to $P$.
\end{enumerate}
No $3$-dimensional platonic solid has all these
properties at the same time.  So, in a sense,
$P$ is even more symmetric than the familiar
$3$-dimensional platonic solids.

The $24$ facets of $P$ are regular octahedra.
These $24$ facets are parallel in pairs, and
there are $12$ codimension $1$ subspaces
such that any facet of $P$ is parallel to
one of these subspaces.  Let $H$ denote
the collection of these subspaces.
The $E_4$-{\it Weyl group\/} is the group 
$W$ generated by reflections in the 
members of $H$.

We make all the same constructions for $P^*$
and we arrive at the group $W^*$ generated 
by reflections in the subspaces parallel to
the facets of $P^*$.  The groups $W$ 
and $W^*$ are conjugate.  We prefer to work
with $W^*$, but we will also consider pay
attention to the subset $H$ of hyperplanes
defined in connection with $W$.

The faces of $P$ have a $3$-coloring such that
parallel faces get the same color such that
each monochrome subset of $H$ consists of
$4$ pairwise perpendicular subspaces.
Concretely, 
we can write $H=H_1 \cup H_2 \cup H_3$
where these set have the following
normal vectors.
\begin{enumerate}
\item $(-a,a,0,0), (0,0,a,-a), (0,0,a,a), (a,a,0,0)$.
\item $(0,a,0,-a), (0,a,0,a), (a,0,-a,0), (a,0,a,0)$.
\item $(a,0,0,a), (-a,0,0,a), (0,a,a,0), (0,a,-a,0)$.
\end{enumerate}

\begin{lemma}
The action of $W^*$ preserves the coloring of the
facets of $P$.
\end{lemma}

\startproof
By symmetry, it suffices to check this for of
the generators of $W^*$.  One of the generators
has the action $(x_1,x_2,x_3,x_4) \to (-x_1,x_2,x_3,x_4)$.
Inspecting the lists above, one can see that this
map preserves the coloring.  Hence, they all do.
\endproof

Let us now turn to the definition of the complex octagonal
PETs.  We will interpret these systems as living in
$\R^4$, using the identification given in
Equation \ref{ident}.  For ease of notation, we
set $X_1=X_1^{\C}$, etc.  Thus $X_1$ and $X_2$ are
$4$-dimensional real polytopes which depend on a
parameter $s$, and $\Lambda_1$ and $\Lambda_2$ are
lattices in $\R^4$ which also depend on $s$.
We will mention $3$ connections with $E_4$.
\newline
\newline
{\bf First Connection:\/}
As we already mentioned, we have $\Lambda_1=\Lambda_2=\Lambda$,
the $E_4$ lattice, when $s=1$.  This is the first connection.
\newline
\newline
{\bf Second Connection:\/}
Let $M^*$ denote the inverse transpose
of the matrix $M$.
The normals to the facets of $X_1$ are the
columns of $(\gamma M_1)^*$, where
$\gamma M_1$ is as in
Equation \ref{marked1}.  The matrix in question is
\begin{equation}
\frac{1}{2}\left[\matrix{-1&0&1/s&1/s\cr 1&0&0&0\cr
0&-1&-1/s&1/s\cr 0&1&0&0}\right]
\end{equation}
Two of these normals appear on the list for $H_1$ and two
of them appear on the list for $H_2$.  Similarly,
two of the normals to the parallelotope
$X_2$ appear on $H_1$ and the other two
appear on $H_2$.
In other words, the $8$ hyperplanes through the
origin parallel to $X_1$ and $X_2$ are
precisely the $8$ hyperplanes in $H_1 \cup H_2$.

Let $A$ be a periodic island for some complex octagonal PET.
By construction, each facet of $A$ is parallel to one of
the facets of $X_1$ or one one of the facets of $X_2$.
Therefore, every face of $A$ is parallel to one of the
hyperplanes in $H_1 \cup H_2$. 
\newline
\newline
{\bf Third Connection:\/}
We say that a polytope $P$ is
$E_4$-{\it semiregular\/} if $P$ is invariant
under the symmetry group $W^*$.
The fixed island for the complex octagonal
PET is precisely $Y=X_1 \cap X_2$.  
Here we prove that $Y$ is $E_4$-semiregular.
This fact establishes the portion of
Conjecture \ref{EE4} having to do with the
fixed point set.

Referring to the discussion in \S \ref{cc},
we see that $Y$ is invariant under the group
of order $8$ generated by the maps
$$(z_1,z_2) \to (\overline z_1,\overline z_2),
\hskip 30 pt
(z_1,z_2) \to (iz_1,iz_2).$$  (We find
it convenient to use complex notation for
the moment.)  Moreover, the map
$$(z_1,z_2) \to (-z_2,z_1)$$ interchanges $X_1$
and $X_2$.  All these maps generate an
order $16$ subgroup $W' \subset W^*$ of
symmetries of $Y$.

Now, $Y$ has at most $16$ sides, and these
sides must come in parallel pairs.  Furthermore,
the sides of $Y$ are colored (say) red
and blue, according as they are parallel
to hyperplanes in $H_1$ or hyperplanes in $H_2$.
The group $W'$ transitively permutes 
the hyperplanes in $H_1$ and also the
hyperplanes in $H_2$.  Therefore,
$W'$ transitively permutes the pairs of
red facets of $Y$ and also transitively
permutes the pairs of blue facets.

But now we can say that all the red facets of $Y$
are the same distance from the origin, and all the
blue facets of $Y$ are the same distance from the
origin.   Since these facets are parallel to the
hyperplanes in $H_1$ and $H_2$, and $W^*$ preserves
each of these sets of hyperplanes, we see that
$W^*$ permutes the hyperplanes extending the red
facets of $Y$.  The same goes for the blue facets.
But $Y$ is just the intersection of halfspaces
bounded by these hyperplanes.  Hence $W^*$ preserves $Y$.
This completes the proof.

\subsection{The Weyl Pseudogroup Action}
\label{E4symm}

After conjecturing that almost every point in a
complex octagonal PET is periodic, one might
be tempted to conjecture that every periodic
island is $E_4$-semiregular.  This, however, is
not the case.  While many of the periodic islands
are $E_4$-semiregular, there are some
anomolous tiles.  We think that the existence
of these anomalous tiles is related to the
red islands discussed in connection with
Figures 1.1 and 1.2.

In the same way that it appears that we can subdivide
the red island in the plane to reveal a more
symmetric pattern for the square turning map, 
we think that perhaps we can
subdivide the asymmetric tiles in the complex
octagonal PET to produce a finer tiling by $E_4$
semiregular polytopes.  If this is really the case,
then one might ask if there is another dynamical
system which produces this finer tiling.
In this section we are going to describe a new
system which seems to produce a finer tiling
and yet seems to be compatible with the complex octagonal
PETs.  This new system is not a mapping in the
traditional sense, but rather a pseudogroup
action.

In this section, we will do three things.  First,
we will describe what we mean by a pseudogroup
action. Second, we will describe the specific
pseudogroup action that is related to the complex
octagonal PETs and give evidence for the connection.
Finally, we will say a word about how we discovered
this alternate system.
\newline
\newline
{\bf Pseudogroup Actions and Orbits:\/}
Let $X \subset \R^d$ be a compact
domain and let
$g_1,...,g_n$ be some finite list of 
isometries.  We insist that this list is
symmetric with respect to inverses.
In other words, an isometry appears
on the list if and only if the inverse
isometry appears on the list.
We do not require that 
$g_k(X)=X$.  

For each $k$, we define
\begin{equation}
X_k=g^{-1}(g_k(X) \cap X).
\end{equation}
By definition, $X_k \subset X$ is the maximal
subset of $X$ such that $g_k(X_k) \subset X$.
We call a subset $Y \subset X$ {\it stable\/}
if it has for following property.  If
$p \in Y \cap X_k$ then $g_k(p) \in Y$ as well.
In other words, the pseudogroup action maps
$Y$ into itself.   We call $Y$ a {\it pseudogroup orbit\/}
if $Y$ is stable and no proper subset of $Y$ is stable.
In case all the isometries preserve $X$, a
pseudogroup orbit corresponds with an orbit
of the group generated by the isometries.
\newline
\newline
{\bf The Weyl Pseudogroup:\/}
Our constructions, as usual, depend on
a parameter $s$.
Consider the matrix
\begin{equation}
A=\frac{s}{2}\left[\matrix{
+1&+1&+1&+1 \cr
+1&+1&-1&-1 \cr
+1&-1&-1&+1 \cr
+1&-1&-1&-1}\right]
\end{equation}
This matrix is a multiple of a Hadamard matrix.
The rows and columns are orthogonal and all
have the same length, namely $s$.

Let $\Pi$ denote the tiling of $\R^4$ by unit
cubes, such that the origin is a vertex of
one of the cubes.  Let $\Gamma_1$ denote the
infinite list of elements generated by reflections
in the hyperplanes of $\Pi$.  By construction,
$\Gamma_1$ is an infinite Bieberbach group.
Let 
\begin{equation}
\Gamma_2=A \circ \Gamma_1 \circ A^{-1}.
\end{equation}
Then $\Gamma_2$ is the group generated by
reflections in the hyperplanes of the
smaller and rotated grid $A(\Pi)$.
No matter how $s$ is chosen, there are only
finitely many elements $g \in \Gamma_1 \cup \Gamma_2$
such that $g(X) \cap X \not = \emptyset$.
To be clear, the group generated by
$\Gamma_1 \cup \Gamma_2$ is a dense group,
and there are infinitely many elements of
this group within every neighborhood of
the identity.  However, we are only allowed
to take elements in $\Gamma_1 \cup \Gamma_2$.

We call $g \in \Gamma_1 \cup \Gamma_2$ {\it good\/}
if either $g(X) \cap X \not = \emptyset$ or
$g^{-1}(X) \cap X \not = \emptyset$.  We call
the pseudogroup generated by the good elements
the {\it Weyl pseudogroup\/}.

To discuss the properties of the Weyl pseudogroup,
we introduce some notation.  Let $O_p$ denote the
orbit of a point $p$ under the complex octagonal PET.
Let $O_p^*$ denote the pseudogroup orbit of $p$.
Let $Q_p$ denote the periodic island about $p$, with
respect to the complex octagonal PET.  Numerically
we observe the following things.

\begin{enumerate}
\item For almost every $p$, the orbit $O_p^*$ is
finite.
\item For almost every $p$, there is a special
point $q \in Q_p$ such that $O_q \subset O_q^*$.
Generically, $O_p^*$ is $192$ times as large as 
$O_q^*$.
\item If $p \in \R^2$ then it seems
$O_q^* \cap \R^2 = O_q$.  
\end{enumerate}
The third item is especially surprising because
$\R^2 \cap X$ is an invariant slice of the
complex octagonal PET -- it is, by definition,
the real octagonal PET.  On the other hand,
$\R^2 \cap X$ is not invariant for the
Weyl pseudogroup.

The properties above might lead the reader to think that the
Weyl pseudogroup is just a minor tweak of the
complex octagonal PET.  However this is not the case.
The pseudogroup orbits tend to be much larger than
the octagonal PET orbits, and indeed the ratio
between the sizes is unbounded.  Nonetheless, the
three properies above suggest that there is an
invariant tiling associated to the Weyl pseudogroup
and that this tiling refines the periodic tiling
associated to the complex octagonal PETs.

The linear parts of the isometies on the
Weyl pseudogroup belong to the group $W^*$
mentioned above.  This property leads us to
believe that the tiles in the invariant tiling
associated to the Weyl pseudogroup are all
$E_4$-semiregular.  We have yet to investigate
this tiling computationally.
\newline
\newline
{\bf The Origins of this System:\/}
Looking at Figures 6.2-6.4, the reader can
probably see that certain subsets of the tiling
exhibit bilateral symmetry.  For instance, the
tilings are all ``as symmetric as possible'' about
$x=k$ for $k=-1,0,1$.  What this means is that
the tiling is symmetric in the domain
$\rho(X) \cap X$, where $\rho$ is a vertical
reflection in one of the lines mentioned.
This is the symmetry we mentioned at the end
of \S \ref{moredetail}.  

We decided to look at the pseudogroup 
generated by all these ``partial symmetries''
of the tiling in the real case.  Eventually
we got the idea of complexifying the picture,
and we arrived at the Weyl pseudogroup
mentioned here.  

\newpage

\section{References}

[{\bf AE\/}] Shigeki Akiyama and Edmund Harris,
{\it Pentagonal Domain Exchange\/}, preprint, 2012
\newline
\newline
[{\bf AG\/}] A. Goetz and G. Poggiaspalla, {\it Rotations by $\pi/7$\/}, Nonlinearity {\bf 17\/}
(2004) no. 5 1787-1802
\newline
\newline
[{\bf AKT\/}] R. Adler, B. Kitchens, and C. Tresser,
{\it Dynamics of non-ergodic piecewise affine maps of the torus\/},
Ergodic Theory Dyn. Syst {\bf 21\/} (2001) no. 4 959-999
\newline
\newline
[{\bf BC\/}] N. Bedaride and J. Cassaigne, {\it Outer Billiards outside regular polygons\/}, J. London Math Soc. (2011) 
\newline
\newline
[{\bf B\/}] J, Buzzi, {\it Piecewise isometries have zero topological entropy (English summary)\/} Ergodic Theory and Dynamical Systems {\bf 21\/}  (2001) no. 5 pp 1371-1377
\newline 
\newline
[{\bf CS\/}] J. Conway and N. Sloane, {\it Lattices and Sphere Packing\/}.
\newline
\newline
[{\bf Go\/}] A. Goetz, {\it Piecewise Isometries -- an emerging area of dynamical systems\/}, preprint.
\newline
\newline
[{\bf GH1\/}] E Gutkin and N. Haydn, {\it Topological entropy of generalized polygon exchanges\/}, Bull. Amer. Math. Soc., {\bf 32\/} (1995) no. 1., pp 50-56
\newline
\newline
[{\bf GH2\/}] E Gutkin and N. Haydn, {\it Topological entropy polygon exchange transformations and polygonal billiards\/}, Ergodic Theory and Dynamical Systems {\bf 17\/} (1997) no. 4., pp 849-867
\newline
\newline
[{\bf GS\/}] E. Gutkin and N. Simanyi, {\it Dual polygonal
billiard and necklace dynamics\/}, Comm. Math. Phys.
{\bf 143\/}:431--450 (1991).
\newline
\newline
[{\bf H\/}] H. Haller, {\it Rectangle Exchange Transformations\/}, Monatsh Math. {\bf 91\/}
(1985) 215-232
\newline
\newline
[{\bf Hoo\/}] W. Patrick Hooper, {\it Renormalization of Polygon Exchage Maps arising from Corner Percolation\/} Invent. Math. 2012.
\newline
\newline
[{\bf K\/}] M. Keane, {\it Non-Ergodic Interval Exchange Transformations\/}, Israel Journal of Math, {\bf 26\/}, 188-96 (1977)
\newline
\newline
[{\bf LKV\/}] J. H. Lowenstein, K. L. Koupsov, F. Vivaldi, {\it Recursive Tiling and Geometry of piecewise rotations by $\pi/7$\/}, nonlinearity {\bf 17\/} (2004) no. 2.
\newline
\newline
[{\bf Low1\/}] J. H. Lowenstein, {\it Aperiodic orbits of piecewise rational
rotations of convex polygons with recursive tiling\/}, Dyn. Syst. {\bf 22\/}
(2007) no. 1 25-63
\newline
\newline
[{\bf Low2\/}] J. H. Lowenstein, {\it Pseudochaotic kicked oscillators\/},
Springer (2012)
\newline
\newline
[{\bf R\/}] G. Rauzy, {\it Exchanges d'intervalles et transformations induites\/},
Acta. Arith. {\bf 34\/} 315-328 (1979)
\newline
\newline
[{\bf S1\/}] R. E. Schwartz, {\it The Octagonal PETs\/}, research
monograph (2012).
\newline
\newline
[{\bf S2\/}] R. E. Schwartz, {\it Outer Billiards on Kites\/}, 
Annals of Math Studies {\bf 171\/} 2009
\newline
\newline
[{\bf S3\/}] R. E. Schwartz, 
{\it Outer Billiards, Quarter Turn Compositions, and Polytope Exchange Transformations\/}, preprint, 2010.
\newline
\newline
[{\bf T\/}] S. Tabachnikov, {\it Billiards\/}, Soci\'{e}t\'{e} Math\'{e}matique de France, 
``Panoramas et Syntheses'' 1, 1995
\newline
\newline
[{\bf W\/}] S. Wolfram, {\it The Mathematica Book\/}, 4th ed., Wolfram Media/Cambridge University Press,
Champaign/Cambridge (1999).
\newline
\newline
[{\bf Y\/}] J.-C. Yoccoz, {\small {\it Continued Fraction Algorithms for Interval Exchange Maps: An Introduction\/} , Frontiers in Number Theory, Physics, and Geometry Vol 1, P. Cartier, B. Julia, P. Moussa, P. Vanhove (ed.) Springer-Verlag 4030437 (2006)\/}
\newline
\newline
[{\bf Z\/}] A. Zorich, {\small {\it Flat Surfaces\/}, Frontiers in Number Theory, Physics, and Geometry Vol 1, P. Cartier, B. Julia, P. Moussa, P. Vanhove (editors) Springer-Verlag 4030437 (2006)\/}

\end{document}